\documentclass[10pt]{article}

\usepackage{psfrag}

\usepackage{graphicx}
\usepackage{epsfig}
\usepackage{amssymb}
\usepackage{amsmath}
\usepackage{amsfonts}
\usepackage{dsfont}
\newcommand{\RR}{\mathbb{R}}

\newcommand{\ZZ}{\mathbb{Z}}
\newcommand{\PP}{\mathbb{P}}
\newcommand{\EE}{\mathbb{E}}

\newcommand{\dps}{\displaystyle}
\newcommand{\eps}{\varepsilon}

\begin{document}

\title{Examples of computational approaches \\ to accommodate randomness in elliptic PDEs}
\author{Claude Le Bris$^\dagger$ \& Fr\'ed\'eric Legoll
\\
Ecole des Ponts and INRIA,
\\
77455 Marne-la-Vall\'ee Cedex 2, France 
\\
{\tt \{claude.le-bris,frederic.legoll\}@enpc.fr}
\\
{\small $\dagger$ Corresponding author}}

\maketitle

\begin{abstract} 
We overview a series of recent works addressing numerical simulations of partial differential equations in the presence of some elements of randomness. The specific equations manipulated are linear elliptic, and arise in the context of multiscale problems, but the purpose is more general. On a set of prototypical situations, we investigate two critical issues present in many settings: variance reduction techniques to obtain sufficiently accurate results at a limited computational cost when solving PDEs with random coefficients, and finite element techniques that are sufficiently flexible to carry over to geometries with random fluctuations. Some elements of theoretical analysis and numerical analysis are briefly mentioned. Numerical experiments, although simple, provide convincing evidence of the efficiency of the approaches. 
\end{abstract}



\section{Introduction}

We consider in this review article a series of works~\cite{banff,mprf,cedya,loz1,loz3,loz2,sqs,minvielle,dcds,minvielle-these}, completed in collaboration with some colleagues of ours, that all share the following common denominator. The task to perform is, possibly repeatedly, to approximate numerically the solution to a partial differential equation that has some random character. In most of our works, the equation has the simplest possible form: it is scalar-valued, elliptic, linear, non degenerate, in divergence form. Typically, and with self-explanatory notations, it reads as 
\begin{equation}
\label{eq:generale}
-\hbox{\rm div} (a \, \nabla u)=f
\end{equation}
on a domain~${\mathcal D}$, and for a certain right hand side $f$. The random character of the problem can be encoded in the coefficient~$a$, and/or in the right-hand side~$f$, and/or in the domain~${\mathcal D}$ itself.

Needless to say, there exist a number of successful approaches to explicitly treat randomness in such partial differential equations. Recent years have witnessed an explosion of the number of methods invented in this extremely lively topic, in particular motivated by the field of \emph{uncertainty quantification}. Stochastic finite elements, spectral methods, sparse tensor products methods, reduced basis techniques, quantization, all methods based one way or another on Karhunen-Loeve, Polynomial Chaos, or other types of similar --or not-- economical decompositions of the random functions present in the equation, both as parameters and unknown functions, have been increasingly studied and considerably improved lately. Some accessible general references in the field are the textbooks~\cite{ghanem,lemaitre}. Review articles on each of the many categories of approaches, such as~\cite{schwab}, are also available. 
The \emph{rationale} behind all these methods is the reduction of the dimension of the problem, considered as a problem in a high-dimensional space consisting of the ambient physical space where the original problem is posed \emph{augmented} by the space of approximation for the random dependency. A simplification of the random dependency follows, and the problem becomes amenable to efficient computational techniques. We ourselves have used, in~\cite{boyo2}, and reviewed, in~\cite{boyo1}, some of these methods (specifically, reduced basis type methods) in the context of equations with random parameters. 

In the present review article, we would like to concentrate on a somewhat alternate strategy: \emph{not attempt to simplify the dependency upon randomness, but embrace the difficulty arising from it}.
Of course, this may only be achieved in some specific situations, sufficiently general to be of broad interest, but certainly not covering the immense spectrum of applications in the engineering and life sciences, and perhaps not with the same generic character as the above mentioned general purpose methods. To some extent, the cases we consider must be slightly simpler than, and not as general as, the cases targeted by the above methods. Our ambition and our achievements are more modest. On the other hand, it turns out that the cases we consider
 arise from a context that brings an additional level of complexity: they originate from multiscale modeling. In that respect, the mix between the presence of randomness, the multiscale feature, and the wish to compute accurately and efficiently lead to an essentially unsolved difficulty (despite the many efforts of outstanding contributors). In multiscale computational science, a number of techniques exist, and are improved constantly. But, most of the times, and despite some overly optimistic claims, they do not marry so well with randomness when it comes to practical computations. Conversely, the extremely efficient toolbox available for random problems has essentially no intersection with multiscale science. Well beyond the somewhat limited purpose of this review article, our intention is thus to attract the attention of the community to the state of the art: \emph{efficiently computing when randomness and many scales are simultaneously present is still an essentially open, considerably challenging, issue}. 

\medskip

\paragraph{The two situations considered: similarities and differences.}

We will consider below two prototypical situations. As we mentioned, we encountered both situations in our research efforts devoted to multiscale science, and more precisely in our endeavor to study and improve computational approaches in materials science: (i) the approximation of the homogenized tensor in the context of stochastic numerical homogenization, and (ii) the multiscale finite element computation of the solution to an (harmless) equation posed on a randomly perforated domain. 

In both situations, randomness originates from geometry. But, in either situation, geometry is encoded differently in the equation. In (i), it is encoded in the heterogeneities of the coefficient. In contrast, in (ii) it is encoded in the computational domain itself. The similarity of the two problematics is evident. By penalization, the second problem may even be viewed as a particular case of the first problem. However, the techniques we employ are different in nature, and have a different purpose. 

\medskip

The former problem (approximation of the homogenized tensor) consists in the repeated resolution of an elliptic equation (the celebrated "corrector equation") on an as large as possible bounded domain, truncation of the whole space $\RR^d$, typically for $d=2$, or $d=3$, in practice. The equation is of the above type~\eqref{eq:generale} and will be made explicit below, see~\eqref{eq:correcteur-random}. The purpose of that repeated solution procedure is to compute an expectation (thus, in practice an empirical mean) from the solutions obtained for various realizations of the local environment (mathematically encoded in the coefficient~$a$ of~\eqref{eq:generale}). The reader may think of different microstructures of the material, different inclusions in the medium, etc. The task of computing an average, that is a single, deterministic output, sounds simple. In particular, the very process of averaging, performed here in the context of a stationary ergodic problem (see the details below), is essentially of the same nature as the law of large numbers. This suggests that the random character of the problem progressively vanishes when the number of realizations considered increases. However, the practical Monte Carlo approach (generate random environments and average out an outcome based on the computed solutions) is plagued by variance issues. The rate of convergence of the approximation, in terms of the size of the truncated computational domain asymptotically covering the whole ambient space, is universal. It is dictated by the central limit theorem. The prefactor appearing in the error estimate is related to the variance of the problem. Efficient computational approaches consist in designing tools to reduce that variance, thus the statistical error in the approximation. That error largely dominates the bias (the deterministic part of the error), and thus is the critical quantity that governs the overall quality of the numerical approach. Our series of works addresses various techniques to reduce the variance: antithetic variables, control variate, selection approach, the latter being somewhat in the spirit of stratified sampling. They are imported from several different contexts and are adjusted to the specific context of numerical homogenization.

\medskip

The latter problem (multiscale finite element computations on randomly perforated domains) differs from the former problem in several respects. As mentioned above, the randomness of the geometry is now encoded not in the coefficients of the equations (which are constant, for simplicity), but in the domain where the equation is posed. We deal with perforations of that domain that are randomly located. Note however that, as briefly mentioned above, the perforations of the domain can be, in the numerical discretization, treated by penalization, in which case the two situations we consider become closer to one another. Another, more substantial, difference is that this second example is a representative case of modern multiscale techniques. In contrast with classical techniques which (i) aim at computing an equivalent "homogenized" medium, and (ii) achieve this in the asymptotic limit of vanishing length scales of the oscillations originally present (this is the case of our former situation), modern techniques attack the multiscale problem (i) directly and (ii) at the actual length scales. Homogenization theory is then seen as a guideline to construct suitable basis functions for the discretization which, in a second stage, are used to expand the numerical solution. Now, in the first stage of construction of these basis functions (the {\it ad hoc} finite elements), the geometry of the computational domain matters. Well before being able to average solutions and/or outputs based on solutions, the issue arises to be able to compute efficiently and accurately each and every single solution, in a robust way that does not require regenerating a set of basis functions for each realization of the random environment. This will be the issue we investigate, and solve using adequate (Crouzeix-Raviart type) boundary conditions for the multiscale finite element basis set and an adequate enrichment (using bubble functions) of that basis set.

We wish to emphasize the following fact. Not only the two problems we consider are relevant beyond the context of multiscale computing, and the methods we put in action can be useful elsewhere, but the second problem considered (robustness of a numerical approach with respect to some wild modifications of the geometry) is not specific to random modeling and covers the generality of all \emph{disordered} geometries, random or deterministic. Not every non periodic problem is indeed random. Techniques literally based upon perfectness of the disorder (such as the approaches we briefly mentioned at the beginning of this introduction, and also the variance reduction techniques we use) can be inefficient in the presence of only imperfect, partial disorder. It is indeed true that the random setting often comes with assumptions such as stationarity, etc, that make randomness tractable at the price of somehow \emph{"rigidifying the disorder"}, if we allow ourselves this abuse of language. Thus the interest of considering other techniques than techniques probabilistic in nature. Other works of ours (see~\cite{m3as} for a review) proceed further in that direction.

\medskip

Before we get to the heart of the matter, let us expose the problems in slightly more details, and announce the plan of our review. 

\paragraph{An elliptic problem with random coefficients.}

The first problem we consider (in Section~\ref{sec:variance} below) arises in the numerical approximation of oscillatory/heterogeneous solutions of PDEs in the context of classical homogenization. The specific setting is that of random ergodic homogenization. The theoretical basis is briefly summarized in Subsection~\ref{ssec:recall}, both in the deterministic (periodic) context and in the random context. The reader willing to familiarize him or herself with that celebrated theory is referred to the classical textbooks~\cite{blp,jikov}, and also to our introductory, computationally oriented lecture notes~\cite{singapour}. The central problem to solve is the corrector problem, which is of the type~\eqref{eq:generale} and more precisely reads as an elliptic equation with a random parameter~$a(x,\omega)$
\begin{equation*}
-\hbox{\rm div}\left[a(x,\omega)\left(p+ \nabla w(x,\omega) \right)\right] = 0
\end{equation*}
(where $p$ is a fixed vector) set on the {\em whole} space $\RR^d$. The precise mathematical formulation is the problem~\eqref{eq:correcteur-random} below. From the solution $ w(x,\omega)$, an expectation is to be evaluated. In practice, the equation is set on a large, truncated domain, subset of $\RR^d$. Several realizations of $a(x,\omega)$ are considered, and an empirical mean is employed to approximate the expectation built from the set of~$ w(x,\omega)$. This is the essence of the Monte Carlo approach in this particular context. Details are given below. Making this Monte Carlo approach efficient requires to reduce the variance. This is the purpose of our series of works~\cite{banff,mprf,cedya,sqs,minvielle,dcds,minvielle-these}, which is overviewed in Subsection~\ref{ssec:variance-reduction}. Well known variance reduction techniques, adapted to the specific context of homogenization, are described and tested. Notice that our research on such topics has already been partially reviewed in~\cite{philtrans,enumath,SMAIProc} where the reader may find more detailed elements of context and some additional comments. In particular, one will also find in~\cite{enumath,SMAIProc} a review of some of our works~\cite{cras-arnaud,arnaud1,arnaud2,jmpa,cras-ronan,m2an-thomines,thomines-redbasis} dealing with problems that are, in a sense made precise there, only "weakly random". Such problems are small random perturbations of periodic problems, and can be approximated by dedicated methods. These "weakly random" problems and methods have a great interest in their own right. They are not addressed in the present review. Our target is "fully random" problems. Some of these methods will however be briefly mentioned below, for different purposes, namely either constructing auxiliary, surrogate models (in Section~\ref{sssec:controlvariate}), or deriving selection rules for the draws (in Section~\ref{sssec:sqs}). Both tools are useful for the reduction of variance in ``fully random'' problems. 

\medskip

In essence, all variance reduction methods aim at eliminating part of the noise intrinsically present in the draws of the realizations of the random variables at play. Noise can be eliminated in a variety of ways. One may try and bias the draws, thereby, in again an informal language, "twisting the arbitrary arm" of randomness. This should be made in a very subtle way. Too small a bias in these draws has no effect on the variance. Too large a bias may indeed succeed in reducing the variance but might affect the average and lead to an incorrect numerical result. Alternately, one may first keep the noise as is, leaving the draws unaffected, but try to cancel the noise out in a second stage, subtracting (in some sense to be made precise) another simulation to the original one. The variance reduction techniques we make use of, which we did not invent but borrowed from other contexts and adapted to ours, proceed in either of the above directions. Notice however that the book of possible variance reduction methods is thick. We have certainly not yet followed all possible tracks. 

The first variance reduction technique we employ is relatively elementary. It is the classical, general purpose technique of \emph{antithetic variables}, presented in Subsection~\ref{sssec:antithetic}. The efficiency of that technique is substantial, but is also, not unexpectedly, limited in particular because the technique does not exploit much the specifics of the problem considered. We present in Subsection~\ref{sssec:controlvariate} the technique of \emph{control variate}, which requires a better knowledge of the problem. A problem simpler to simulate and close to the original problem, in a sense that is made precise below, has to be considered and concurrently solved. The technique uses that knowledge to get a much better reduction of the variance. Specifically, the auxiliary problem used as control is built upon an approximation method that we developed for the "weakly random" problems we were alluding to above. In some sense, this problem is a reduction of the original problem. In Subsection~\ref{sssec:sqs}, we expose yet another approach, specifically imported from solid state physics, namely that of \emph{special quasi-random structures}. It consists in selecting, in the spirit of another well-known technique, {\em stratified sampling}, some configurations of the random environment that are more suitable than generic configurations to compute the empirical averages, so as to again reduce the variance. Similarly to the previous approach of control variate, a reduced model is also useful in the approach.

\medskip

\paragraph{An elliptic problem on a randomly perforated domain.}

Our second framework concerns a PDE set on a perforated domain. Randomness manifests itself in the problem via the randomness of the location of the perforations (and, possibly, their size and shape). Again, we consider the simplest possible equation: scalar-valued, elliptic, non degenerate, linear and in divergence form. More precisely, we consider a bounded domain~${\mathcal D} \subset \RR^d$ and a set $B_\eps$ of perforations within this domain. The perforations are supposedly small and in extremely large a number. The parameter~$\eps$ encodes the typical small distance between the perforations. We denote by~${\mathcal D}_\eps={\mathcal D} \setminus \overline{B_\eps}$ the perforated domain (see Figure~\ref{fig:perforation} in Section~\ref{sec:msfem}). We then consider the following problem: find $u : {\mathcal D}_\eps \rightarrow \mathbb{R}$, solution to
\begin{equation}
\label{eq:genP}
-\Delta u = f \ \text{ in ${\mathcal D}_\eps$},
\quad
u = 0 \ \text{ on $\partial {\mathcal D}_\eps$},
\end{equation}
where $f:{\mathcal D} \rightarrow \mathbb{R}$ is a given function, assumed sufficiently regular on ${\mathcal D}$. It is important to note that the homogeneous Dirichlet boundary condition on $\partial {\mathcal D}_\eps$ (and hence on the boundary ${\mathcal D} \cap \partial B_\eps$ of the perforations) is a crucial feature of the problem we consider. The consideration of this particular academic problem stems from our interest in various physically relevant problems, for instance in fluid mechanics, atmospheric modeling, electrostatic devices, etc. Equation~\eqref{eq:genP} can also be seen as a step toward the resolution of advection-diffusion equations (see our ongoing work~\cite{madiot-these} for that particular category of equations), that of the Stokes problem (see e.g.~\cite{lozinski_stokes}), or that of the Navier-Stokes problem, all on perforated domains. In fluid dynamics applications, homogeneous Dirichlet boundary conditions on the perforations are typical.

\medskip

We introduce and study a dedicated multiscale finite element method (MsFEM). The variant of MsFEM we consider uses \emph{Crouzeix-Raviart type} finite elements~\cite{Crouzeix-Raviart}. We have introduced this approach in~\cite{loz1,loz3} on a multiscale elliptic problem (exactly of the form~\eqref{eq:generale}) with a highly oscillatory coefficient (the coefficient~$a$ in~\eqref{eq:generale}), and then extended it in~\cite{loz2} to the case of problem~\eqref{eq:genP}. In the presence of perforations, we indeed have to improve the previously introduced approach by the addition of \emph{bubble functions} to the finite element basis set. Both ingredients, Crouzeix-Raviart type finite elements on the one hand, and addition of bubble functions on the other hand, allow for robustness of the approach, a critical issue when the geometry of the computational domain is random.

It is well known that, for the construction of multiscale finite elements, boundary conditions set on the edges (or facets in 3D) of mesh elements for the definition of the basis functions play a critical role for the eventual accuracy and efficiency of the approach. Using Crouzeix-Raviart type elements (see~\cite{Crouzeix-Raviart} for their original introduction in the context of classical --meaning non-multiscale-- finite elements) gives a definite flexibility. In short, the continuity of the finite element basis set functions across the edges of the mesh is enforced only in a weak sense by requiring that the average of the jump vanishes on each edge. The nonconforming approximation obtained in this manner proves to be very effective. The above issue regarding boundary conditions on the mesh elements is all the more crucial when dealing with perforated computational domains, and may be of paramount importance for \emph{randomly} perforated computational domains. The approach is meant to be as insensitive as possible to the possible intersections between element edges and perforations. The perforations can then be very heterogeneously distributed. An {\it ad hoc} construction of a mesh that as much as possible avoids intersecting the perforations is indeed prohibitively difficult. Our approach does {\em not} require a mesh of this type. 

The second ingredient of our approach is the addition of bubble functions to the finite element space. As illustrated using a simple one-dimensional analysis (see~\cite{loz2}), and demonstrated with the extensive set of numerical tests performed in~\cite{loz2}, some of which will be reproduced here, bubble functions definitely bring an added value for the overall accuracy of the approach. 

\medskip

In the vast literature devoted to similar types of problems and techniques, we refer the reader unfamiliar with the general context to~\cite{Brenner-Scott,Quarteroni} for some background on finite element methods in general, to the review article~\cite{arbogast-review} and the book~\cite{Efendiev-Hou-book} (and references therein) for the general background on MsFEM, and the works~\cite{henning,hornung,lions1980} specifically addressing problems on perforated domains, either from a theoretical or a numerical standpoint. 

\medskip

Our contribution is reviewed in Section~\ref{sec:msfem}. That section begins with some generalities about our MsFEM approach, contained in Subsection~\ref{ssec:presentation2D}. In particular, we briefly mention an error estimate that we have established in the case of periodic perforations. We are unfortunately unable to extend it to the case of random perforations. It is however useful to illustrate the nice properties of the approach. Subsection~\ref{ssec:Numerical-tests} is then devoted to some numerical tests that show the announced robustness of the approach in the presence of random perforations. All the numerical tests we report below have been performed using the finite elements software FreeFem++~\cite{hecht}.

\section{Variance reduction for elliptic PDEs with random coefficients}
\label{sec:variance}

\subsection{Our specific context: numerical homogenization}
\label{ssec:recall}

We first recall some basic ingredients of elliptic homogenization theory in the periodic setting. We consider, in a regular domain ${\cal D}$ in $\RR^d$, the problem 
\begin{equation}
\label{eq:pb0-per}
\left\{
\begin{array}{l}
\dps -\hbox{\rm div}\left[ A_{\rm per}\left(\frac{x}{\eps}\right) \nabla u^\eps \right] = f \quad \text{in ${\cal D}$}, 
\\ \noalign{\vskip 5pt}
u^\eps=0 \quad \text{on $\partial {\cal D}$}, 
\end{array}
\right.
\end{equation}
where the matrix $A_{\rm per}$ is symmetric and $\ZZ^d$-periodic. We manipulate for simplicity \emph{symmetric} matrices, but our discussion carries over to non symmetric matrices up to slight modifications.

The corrector problem associated to (\ref{eq:pb0-per}) reads, for $p$ fixed in $\RR^d$,
\begin{equation}
\label{eq:cor-per-intro}
\left\{
\begin{array}{l}
-\hbox{\rm div}\left(A_{\rm per}(y)\left(p+ \nabla w_p(y) \right)\right) = 0, 
\\ \noalign{\vskip 5pt}
w_p \text{ is $\ZZ^d$-periodic}.
\end{array}
\right.
\end{equation}
It has a unique solution up to the addition of a constant. Then, the homogenized coefficients read $\dps [A^\star]_{ij} = \int_Q e_i^T A_{\rm per}(y) \left(e_j + \nabla w_{e_j}(y) \right)dy$, where $Q$ is the unit cube and $\left\{ e_i \right\}_{1 \leq i \leq d}$ is the canonical basis of $\RR^d$. The main result of periodic homogenization theory is that, as $\eps$ goes to zero, the solution $u^\eps$ to (\ref{eq:pb0-per}) converges to $u^\star$ solution to
\begin{equation}
\label{eq:pb0-star}
\left\{
\begin{array}{l}
-\hbox{\rm div}\left[ A^\star \nabla u^\star \right] = f \quad \text{in ${\cal D}$}, 
\\ \noalign{\vskip 5pt}
u^\star=0 \quad \text{on $\partial {\cal D}$}.
\end{array}\right.
\end{equation}
The convergence holds in $L^2({\cal D})$, and weakly in $H^1_0({\cal D})$. The correctors $w_{e_i}$, $1 \leq i \leq d$, may then also be used to ``correct'' $u^\star$ in order to identify the behavior of $u^\eps$ in the strong topology $H^1_0({\cal D})$. Practically, at the price of only computing $d$ periodic problems (\ref{eq:cor-per-intro}), the solution to problem (\ref{eq:pb0-per}) can be efficiently approached for $\eps$ small.

\medskip

The random setting is, for the simple equation~\eqref{eq:pb0-per}, a straightforward extension of the periodic setting. Skipping all technicalities related to the definition of the probabilistic setting, which we assume discrete stationary and ergodic (we refer e.g. to~\cite{singapour} for all details), we now fix $A(\cdot,\omega)$ a square matrix of size $d$, which is assumed stationary in the sense that
\begin{equation}
\label{eq:stationnarite-disc}
\forall k\in \ZZ^d, \quad A(x+k, \omega) = A(x,\tau_k\omega) \mbox{ almost everywhere in $x$, almost surely}
\end{equation}
(where~$\tau$ is an ergodic group action) and which is assumed to enjoy the classical assumptions of uniform ellipticity and boundedness. Then we consider the boundary value problem
\begin{equation}
\label{eq:pb0-stoc}
\left\{
\begin{array}{l}
\dps -\hbox{\rm div}\left(A\left(\frac{x}{\eps}, \omega\right) \nabla u^\eps \right) = f \quad \text{in ${\mathcal D}$}, 
\\ \noalign{\vskip 5pt}
u^\eps = 0 \quad \text{on $\partial {\mathcal D}$}.
\end{array}
\right.
\end{equation}
Standard results of stochastic homogenization~\cite{blp,jikov} apply and allow to find the homogenized limit of problem (\ref{eq:pb0-stoc}). The solution $u^\eps(\cdot,\omega)$ to (\ref{eq:pb0-stoc}) converges to the solution to (\ref{eq:pb0-star}) where the homogenized matrix is now defined as
\begin{equation}
\label{eq:Astar-random}
[A^\star]_{ij} 
= 
\EE\left(\int_Q e_i^T A(y,\cdot) \, \left( e_j+\nabla w_{e_j}(y,\cdot)\right) \,dy\right),
\end{equation}
where $Q$ is the unit cube and where, for any $p\in \RR^d$, $w_p$ is the solution (unique up to the addition of a random constant) to
\begin{equation}
\label{eq:correcteur-random}
\left\{
\begin{array}{l}
-\hbox{\rm div}\left[A(y,\omega)\left(p+ \nabla w_p(y,\omega) \right) \right] = 0 \quad \hbox{\rm a.s. on $\RR^d$}, 
\\ \\
\nabla w_p \quad \mbox{is stationary in the sense of (\ref{eq:stationnarite-disc})}, 
\\ \\
\dps \EE\left(\int_Q \nabla w_p(y,\cdot)\,dy\right) = 0. 
\end{array}
\right.
\end{equation}
A striking difference between the stochastic setting and the periodic setting can be observed comparing (\ref{eq:cor-per-intro}) and (\ref{eq:correcteur-random}). In the periodic case, the corrector problem is posed on a bounded domain (namely, the periodic cell~$Q$), since the corrector $w_p$ is periodic. In sharp contrast, the corrector problem (\ref{eq:correcteur-random}) of the random case is posed on the whole space $\RR^d$. In the numerical practice, truncations of that problem have to be considered, typically on large domains~$Q_N$ (say, $Q_N = (-N,N)^d$) and using periodic boundary conditions:
\begin{equation}
\label{eq:correcteur-random-N}
\left\{
\begin{array}{l}
-\hbox{\rm div}\left[A(y,\omega)\left(p+ \nabla w_p^N(y,\omega) \right)\right] =0,
\\ \\
w_p^N(\cdot,\omega) \mbox{ is $Q_N$-periodic}.
\end{array}
\right.
\end{equation}
The random matrix $A^\star_N(\omega)$ defined by
\begin{equation}
\label{eq:Astar-random-N}
[A^\star_N(\omega)]_{ij} 
= 
\frac{1}{|Q_N|} \int_{Q_N} e_i^T A(y,\omega) \, \left( e_j+\nabla w_{e_j}^N(y,\omega)\right) \,dy
\end{equation}
converges (almost surely) to $A^\star$ when $N \to \infty$ (see e.g.~\cite{bourgeat2}). Important theoretical questions about the quality and the rate of the convergence in terms of the truncation size arise (see, in particular, the pioneering works by A.~Bourgeat and A.~Piatnitski~\cite{bourgeat1,bourgeat2} and a recent series of works by A.~Gloria, F.~Otto and their collaborators~\cite{gloria-otto}). The error
\begin{equation}
\label{eq:error-decomposition}
A^\star - A_N^\star(\omega)
=
\Big( A^\star - \EE \left[ A_N^\star \right] \Big)
+
\Big( \EE \left[ A_N^\star \right] - A_N^\star(\omega) \Big)
\end{equation}
is the sum of a systematic error and of a statistical error (the first and second terms in the above right-hand side, respectively). It is observed, and theoretically demonstrated in some cases, that the statistical error decays with a slower rate with respect to $N$ than the systematic error. For large values of $N$, the statistical error is therefore dominating the systematic error. In what follows, we focus on reducing the statistical error, and describe approaches to more efficiently compute $\EE \left[ A_N^\star \right]$, for a given truncated domain $Q_N$. 

\subsection{Variance issues and their reduction}
\label{ssec:variance-reduction}

The direct approach for numerical random homogenization consists in using empirical means approximating the expectation $\mathbb{E} \left(A^\star_N \right)$ of the homogenized coefficient computed on a truncated domain of size~$N$. Variance issues are central.

\subsubsection{Antithetic variables}
\label{sssec:antithetic}

\paragraph{Description of the approach.} 

\begin{figure}[htbp]
\centering\includegraphics[width=3.5in]{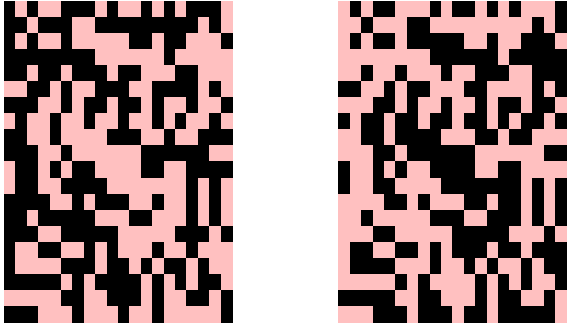}
\caption{An example of realization of $A(x,\omega)$ together with its antithetic field $B(x,\omega)$ for the famous random checkerboard (reproduced from~\cite{cedya}). \label{fig_B12}}
\end{figure}

The variance reduction technique known as the technique of \emph{antithetic variables} has been adapted to the homogenization context in~\cite{banff,mprf,cedya,dcds}. It consists in concurrently considering two sets of configurations for the random material instead of only one set. The two sets of configurations are deduced one from the other (see Figure~\ref{fig_B12} for a simple example). Fix $M=2 \mathcal{M}$ and consider $\mathcal{M}$ i.i.d. copies $\left(A^m(x,\omega)\right)_{1 \leq m \leq \mathcal{M}}$ of $A(x,\omega) $. Construct next $\mathcal{M}$ i.i.d. \emph{antithetic} random fields $B^m(x,\omega)= T \left( A^m(x,\omega) \right)$, $1 \leq m \leq \mathcal{M}$, from the $\left(A^m(x,\omega)\right)_{1 \leq m \leq \mathcal{M}}$. The map $T$ transforms the random field $A^m$ into another, so-called \emph{antithetic}, field $B^m$. The transformation is performed in such a way that, for each $m$, $B^m$ has the same law as $A^m$, namely the law of the matrix~$A$. Somewhat vaguely stated, if $A$ was obtained in a coin tossing game (using a fair coin), then $B^m$ would be {\it head} each time $A^m$ is {\it tail} and vice versa. Then, for each $1 \leq m \leq \mathcal{M}$, we solve two corrector problems of the type~\eqref{eq:correcteur-random-N}. One is associated to the original $A^m$, the other one is associated to the antithetic field $B^m$. Using its solution $v_p^{N,m}$, we define the \emph{antithetic homogenized matrix} $B^{\star,m}_N(\omega)$, the elements of which read, for $1 \leq i,j \leq d$, 
$$
\left[ B^{\star,m}_N(\omega) \right]_{ij}
= 
\frac{1}{|Q_N|} \int_{Q_N} e_i^T B^m(y,\omega) \ \left(e_j + \nabla v_{e_j}^{N,m}(y,\omega) \right) \, dy.
$$
Next we set, for any $1 \leq m \leq \mathcal{M}$,
\begin{equation*}
\widetilde{A}^{\star,m}_N(\omega) 
= 
\frac{1}{2} \left( A^{\star,m}_N(\omega) + B^{\star,m}_N(\omega) \right).
\end{equation*}
Since $A^m$ and $B^m$ are identically distributed, so are $A^{\star,m}_N$ and $B^{\star,m}_N$. Thus, $\widetilde{A}^{\star,m}_N$ is unbiased (that is, $\mathbb{E}\left(\widetilde{A}^{\star,m}_N\right)=\mathbb{E}\left({A}^{\star,m}_N\right)$). In addition, it satisfies $\widetilde{A}^{\star,m}_N \underset{N \rightarrow + \infty}{\longrightarrow} A^\star$ almost surely, because $B^m$ is ergodic. The hope (which is indeed confirmed by theory in simple cases and by numerical observation in all the tests we performed) is that the new approximation~$\widetilde{A}^{\star,m}_N$ has a smaller variance than~$A^{\star,m}_N$. It is indeed the case under appropriate assumptions.

\paragraph{Foundation of the approach.} 

Some theoretical elements that give a foundation to our approach of variance reduction are collected in~\cite{mprf}. Of course, as always in this context, the one-dimensional setting, being explicit, may be entirely analyzed. Two particular higher dimensional cases were also analyzed. The first one is a ``genuinely'' random setting, the second one is a ``weakly random'' case. In all such cases, under appropriate assumptions, the approach is shown to qualitatively reduce the variance. A quantitative assessment of the reduction is however out of reach. Only numerical tests can provide some information.

\paragraph{Numerical tests.} 

\begin{figure}[htbp]
\psfrag{Number of cells: \(2N\)^2}{Size of $Q_N$}
\centering\includegraphics[width=3.5in]{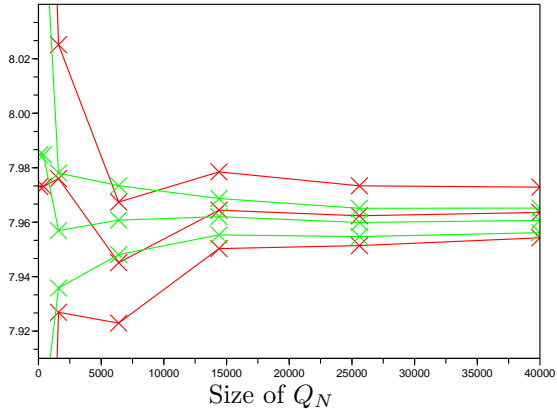}
\caption{Estimation of $A_{11}^\star$ (with confidence interval) with respect to $|Q_N|$ (in red, the classical MC strategy, in green the antithetic variable strategy, for an equal computational cost; reproduced from~\cite{cedya}). \label{fig_B12_result}}
\end{figure}

The tests performed concern three different ``input'' random fields $A(\cdot,\omega)$ in~\eqref{eq:pb0-stoc}, some i.i.d., some correlated, with various correlation lengths. We have investigated variance reduction on a typical diagonal~$\left[A^\star_N \right]_{11}$, or off-diagonal~$\left[A^\star_N \right]_{12}$ entry of the approximate homogenized matrix~$A^\star_N$, as well as on the eigenvalues of the matrix, and the eigenvalues of the associated differential operator~$L_A = - \hbox{\rm div} \left[ A^\star_N(\omega) \nabla \cdot \right]$ (supplied with homogeneous Dirichlet boundary conditions on~$\partial{\cal D}$). We show a typical result on Figure~\ref{fig_B12_result}, for a two-dimensional test case where $A(x,\omega)$ is given by
$$
A(x,\omega) = \sum_{k\in\ZZ^2} {\bf 1}_{Q+k}(x) \, a_k(\omega) \, \text{Id},
$$
where $Q=(0,1)^2$, $\text{Id}$ is the identity matrix, and $\left(a_k\right)_{k\in\ZZ^2}$ is an i.i.d sequence of random variables such that $\PP(a_k = \alpha) = \PP(a_k = \beta) = 1/2$, with $\alpha = 3$ and $\beta = 20$. We observe on Figure~\ref{fig_B12_result} that the confidence interval obtained with our approach is smaller than that obtained with the standard Monte Carlo (MC) approach. The gain in variance (at equal computational cost) is essentially insensitive to the size of the computational domain. It is here between 6 and 10. At fixed computational cost, the approach improves the accuracy by a factor $\sqrt 6 \approx 2.45$. Equivalently, for a given accuracy, the computational cost is reduced by a factor $6$.

\medskip

The approach has been extended in~\cite{dcds} to nonlinear convex stochastic homogenization problems. The highly oscillatory problem, the Euler-Lagrange equation of which replaces~\eqref{eq:pb0-stoc}, is then given by
$$
\inf \left\{ \int_{\cal D} W\left(\frac{x}{\eps}, \omega, \nabla u(x) \right) \, dx - \int_{\cal D} f(x) u(x) dx, \quad u \in W^{1,q}_0({\cal D}) \right\}
$$
for some appropriate $q \geq 2$. We assume that $W$ is strictly convex with respect to its third argument, and that, for any $\xi \in \RR^d$, the random field $(y,\omega) \mapsto W(y,\omega,\xi)$ is stationary in the sense of~\eqref{eq:stationnarite-disc}. We also assume that $W$ satisfies a standard growth condition (see~\cite{dcds} for details). The associated homogenized problem, the Euler-Lagrange equation of which replaces~\eqref{eq:pb0-star}, is
$$
\inf \left\{ \int_{\cal D} W^\star\left(\nabla u(x) \right) \, dx - \int_{\cal D} f(x) u(x) dx, \quad u \in W^{1,q}_0({\cal D}) \right\}
$$
where, for any $p \in \RR^d$, the homogenized energy density $W^\star$ is given by $\dps W^\star(p) = \lim_{N \to \infty} W^\star_N(\omega,p)$ where
$$
W^\star_N(\omega,p) 
= 
\inf \left\{ \frac{1}{|Q_N|} \int_{Q_N} W\left(y,\omega, p + \nabla w(y) \right) \, dy, \quad w \in W^{1,q}_{\rm per}(Q_N) \right\}.
$$
As in the linear case for the corrector problems~\eqref{eq:correcteur-random} and~\eqref{eq:correcteur-random-N}, we approximate in practice $W^\star(p)$ by $\EE \left[ W^\star_N(\cdot,p) \right]$ for some large $N$, the latter quantity being itself approximated by an empirical mean. We show in~\cite{dcds} that the exact same antithetic variable technique we used in the linear case can be used to compute $\EE \left[ W^\star_N(\cdot,p) \right]$ (for any $N$ and any $p \in \RR^d$) more efficiently. 

\bigskip

To summarize, our numerical results show that the technique may be applied to a large variety of situations and has proved efficient whatever the output considered. Variance is systematically reduced. We have observed, however, that the rate of reduction is not spectacular. This has motivated the consideration of other approaches.

\subsubsection{Control variate technique}
\label{sssec:controlvariate}

The \emph{control variate approach} is a variance reduction technique known to be potentially much more efficient than the antithetic variable technique. It however asks to have beforehand a better information on the random variables simulated. In the context of homogenization, the work~\cite{minvielle} presents a first possible investigation of the efficiency of this technique. As briefly mentioned in the introduction, the approach aims at "subtracting" part of the noise, using an auxiliary problem close to the original problem. We describe how we derive such an auxiliary problem in the next paragraph. Next, we make the setting specific, and present the variance reduction approach.

\medskip

\paragraph{A useful ``weakly random'' setting.}

\begin{figure}
\centering
\hskip5truemm\psfig{figure=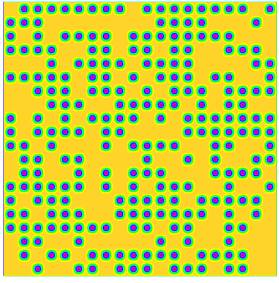,height=4truecm}
\caption{A typical realization of the Bernoulli law for the perturbed periodic material. \label{fig:arnaud-bernoumat}}
\end{figure}

One perturbative approach among many other possible ones, described in full details in~\cite{cras-arnaud,arnaud1,arnaud2}, and addressing the random material as a small perturbation of a periodic material, consists in considering a coefficient of the form
\begin{equation}
\label{eq:perturb2}
A_\eta(x, \omega) = A_{\rm per}(x) + b_\eta(x, \omega) C_{\rm per}(x),
\end{equation}
where, with evident notation, $A_{\rm per}$ is a periodic matrix modeling the unperturbed material, and $C_{\rm per}$ is a periodic matrix modeling the perturbation. Consider then the case
\begin{equation}
\label{consmail}
b_\eta(x,\omega) = \sum_{k \in \mathbb{Z}^d} \mathbf{1}_{Q+ k}(x) \, B_\eta^k(\omega),
\end{equation}
where the $B_\eta^k$ are, say, independent identically distributed scalar-valued random variables. One particularly interesting case (see~\cite{cras-arnaud,arnaud1,arnaud2} for other cases) is when the common law of the $B_\eta^k$ is assumed to be a Bernoulli law of parameter~$\eta$ (see Figure~\ref{fig:arnaud-bernoumat}). We now \emph{formally} explain how we derive an expansion of the problem in ``powers'' of $\eta$. The approach has been recently proven to be rigorous in~\cite{mourrat,gloria-arma}.

\medskip

The basic observation consists in noticing that, in the corrector problem
\begin{equation}
\label{eq:correc-intuition}
-\hbox{\rm div}\left[A_\eta(y,\omega)\left(p+ \nabla w_p(y,\omega) \right)\right] = 0,
\end{equation}
the only source of randomness comes from the coefficient $A_\eta(y,\omega)$. Therefore, in principle, if one knows the law of this coefficient, one knows the law of the corrector function $w_p(y,\cdot)$ and therefore may compute the homogenized coefficient~$A_\eta^\star$, the latter being a (quite intricate but well defined) function of this law. When the law of $A_\eta$ is an expansion in terms of a small coefficient, so is the law of $w_p$. Consequently $A_\eta^\star$ can be obtained as an expansion. Heuristically, on the cube~$Q_N = [0,N]^d$ and at order 1 in $\eta$, the probability to get the perfect periodic material (entirely modeled by the matrix~$A_{\rm per}$) is $(1-\eta)^{N^d}\approx 1-N^d\eta+O(\eta^2)$, while the probability to obtain the unperturbed material on all cells except one (where the material has matrix $A_{\rm per} + C_{\rm per}$) is $\eta \, N^d \, (1-\eta)^{N^d-1} \approx \eta \, N^d +O(\eta^2)$. All other configurations, with more than two cells perturbed, yield contributions of orders higher than or equal to $\eta^2$. This gives the intuition (confirmed by the mathematical proof) that the first order correction indeed comes from the difference between the material perfectly periodic except on one cell and the perfect material itself. The homogenized matrix reads as $A_\eta^\star = A_{\rm per}^\star+ \eta A_{1,\star} + o(\eta)$ where $A_{\rm per}^\star$ is the homogenized matrix for the unperturbed periodic material and
\begin{equation*} 
A_{1,\star}\, e_i = \lim_{N \rightarrow + \infty} \int_{Q_N}\left[(A_{\rm per}+\mathbf{1}_{Q}C_{\rm per})(e_i + \nabla w_i^N) - A_{\rm per}(e_i + \nabla w_i^0) \right], 
\end{equation*}
where $w_i^0$ is the $Q$-periodic corrector for $A_{\rm per}$ (solution to~\eqref{eq:cor-per-intro}), and $w_i^N$ solves
\begin{equation} 
\label{eq:defcell}
-\mathrm{div}\left[ (A_{\rm per}+ \mathbf{1}_{Q}C_{\rm per}) (e_i + \nabla w_i^N) \right] = 0 \quad \text{in $Q_N$}, \quad \text{$w_i^N$ is $Q_N$-periodic}.
\end{equation}

\medskip

The approach has been extensively tested. It is observed that, using the perturbative approach, the large $N$ limit for cubes of size $N$ is already very well approached for small values of $N$. The computational efficiency of the approach is clear: solving the two periodic problems with coefficients $A_{\rm per}$ and $A_{\rm per}+ \mathbf{1}_{Q}C_{\rm per}$ for a limited size $N$ is much less expensive than solving the original, random corrector problem for a much larger size $N$. When the second order term is needed (and we will actually use it below), configurations with two defects have to be computed. They all can be seen as a family of PDEs, parameterized by the geometrical location of the defects. Reduced basis techniques have been shown to be useful and allow for a definite speed-up in the computation, see~\cite{thomines-redbasis}. 

\medskip

\paragraph{Variance reduction in the ``fully random'' setting.}

When the parameter~$\eta$ of the Bernoulli law is \emph{not} taken small in the definition~\eqref{eq:perturb2}-\eqref{consmail} of the coefficient inserted in the equation~\eqref{eq:pb0-stoc}, and thus the corrector problem~\eqref{eq:correcteur-random}, the expansion technique we have just described is not correct, and actually we observe experimentally that it is inaccurate. It can however serve for the construction of a control variate, useful to reduce the variance. 

\medskip

We recall that the field $b_\eta$ writes as in~\eqref{consmail}. Determining the field $A(x,\omega)$, given by~\eqref{eq:perturb2}, on the truncated domain~$Q_N$ amounts to drawing $B^k_\eta(\omega)$ in each cell $Q+k$ in~$Q_N$. This allows to compute the associated (approximate) homogenized coefficient $A^\star_N(\omega)$ from the solution to the corrector problem~\eqref{eq:correcteur-random-N} (with $A \equiv A_\eta$), see~\eqref{eq:Astar-random-N}. In parallel to this task, we \emph{reconstruct} from the specific realization of the set of $\left\{ B_\eta^k(\omega) \right\}_{k \, \text{s.t.} \, Q+k \subset Q_N}$ a field that is used as a control variate. More precisely, we set
\begin{equation}
\label{eq:encore}
C^\star_N(\omega) =A^\star_N(\omega)-\rho \left(A_{\rm per}^\star + A_1^{\star,N}(\omega)-\EE \left[ A_{\rm per}^\star + A_1^{\star,N}(\omega) \right] \right).
\end{equation}
In this formula, 
$$
A_1^{\star,N}(\omega)= \frac{1}{|Q_N|} \sum_{Q+k\subset Q_N} {B^k_\eta(\omega) } \ {\cal A}^{\rm 1 \, def}_k 
$$
where ${\cal A}^{\rm 1 \, def}_k$ is the deterministic coefficient corresponding to the case of \emph{one} defect located at position~$k$ in $Q_N$ (it is actually independent of~$k$ thanks to the periodic boundary conditions in~\eqref{eq:defcell}). The parameter~$\rho$ is a deterministic parameter, a traditional ingredient of control variate techniques, which is optimized in terms of the estimated variances of the objects at play. It is crucial to note that the expectation of~$A_1^{\star,N}(\omega)$ is analytically computable. Since by construction $\EE\left(C^\star_N\right) = \EE\left(A^\star_N\right)$, the technique then consists in approximating the former (thus the latter) by an empirical mean. The theoretical study and the numerical tests in~\cite{minvielle} show that the variance of $C^\star_N(\omega)$ is smaller than that of $A^\star_N(\omega)$, and hence that the quality of the approximation is improved. 

\medskip

\begin{figure}[htbp]
\psfrag{Supercell sidelength N}{\footnotesize \hspace{7mm} $N$}
\psfrag{Homogenized coefficient}{\footnotesize \hspace{-7mm} Homogenized coefficient}
\centering
\includegraphics[width=2.9in]{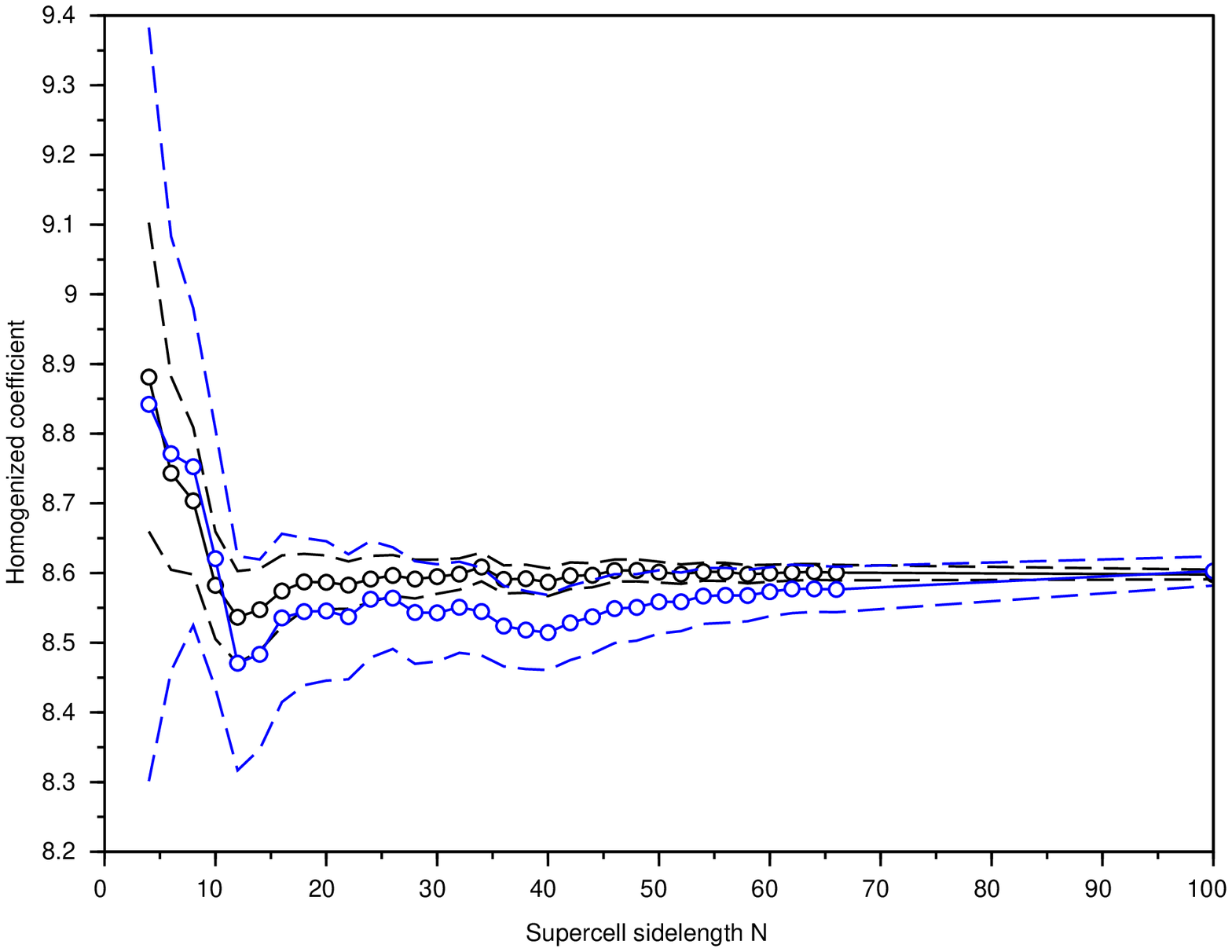} 
\includegraphics[width=2.9in]{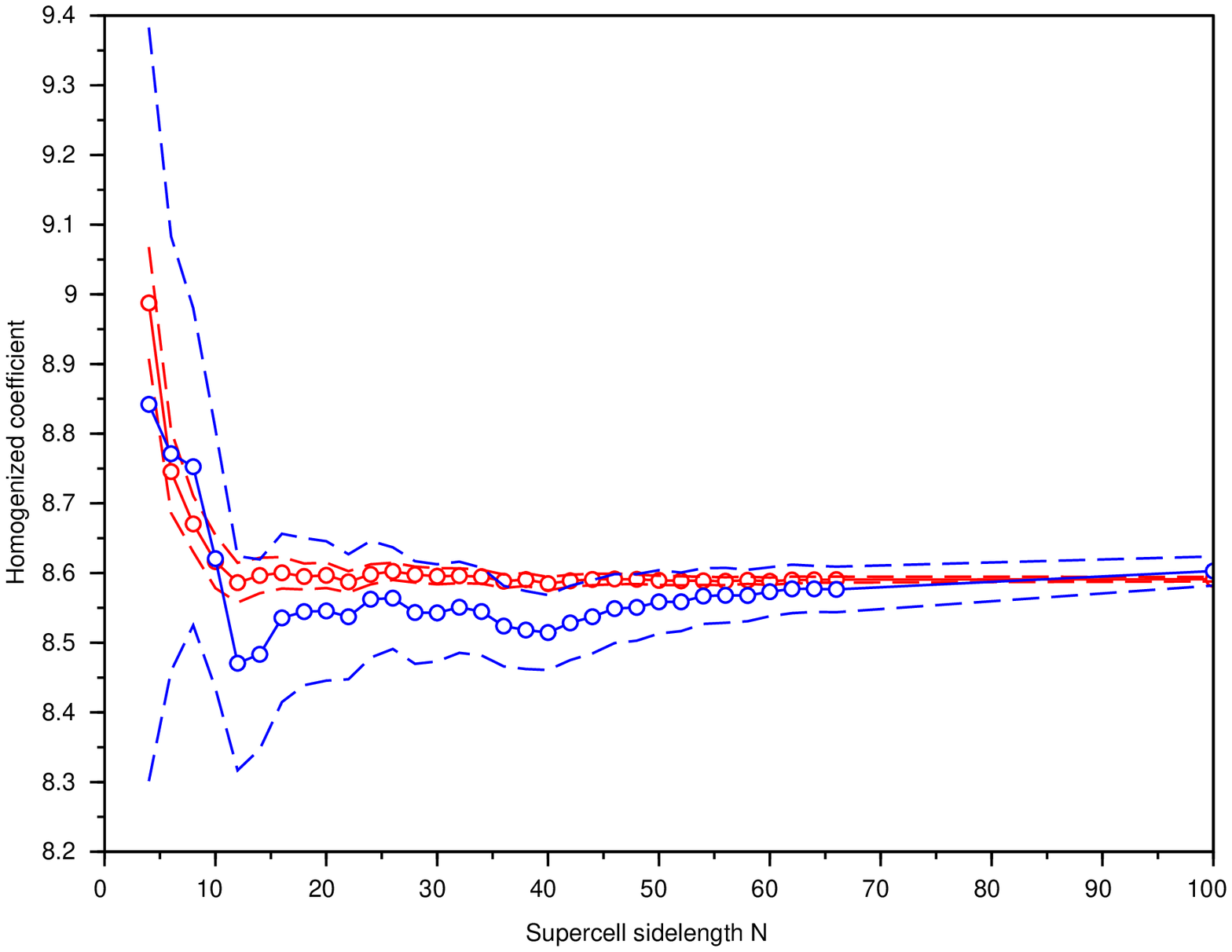}
\caption{Estimation of $A_{11}^\star$ together with its confidence interval (computed using $M=100$ i.i.d. realizations). Left: we use the classical MC approach (in blue) and the control variate approach based on~\eqref{eq:encore} (in black). Right: we use the classical MC approach (in blue) and the second-order control variate approach (in red) (reproduced from~\cite{minvielle}). \label{CV_result}}
\end{figure}

As an illustration, we use a similar case as in Subsection~\ref{sssec:antithetic}, which reads as~\eqref{eq:perturb2}--\eqref{consmail} with $\eta = 1/2$. Applying the above strategy provides the results on the left part of Figure~\ref{CV_result}, where the variance is reduced by a factor close to 6, that is, comparable to the technique of antithetic variables. One may next improve on the previous results using a \emph{second order} expansion with respect to $\eta$ and including in the control variate the deterministic coefficients corresponding to the case of both one and two defects. Additional parameters playing the role of $\rho$ above are needed, in order to ensure a substantial variance reduction (see the details in~\cite{minvielle}). The variance reduction of such a case, of the order of 40, is represented on the right part of Figure~\ref{CV_result}.

\subsubsection{Special Quasi-Random Structures}
\label{sssec:sqs}

The variance reduction approach we now overview has been originally introduced by other authors for a slightly different purpose in atomistic solid-state science~\cite{vonPezoldDickFriakNeugebauer2010,WeiFerreiraBernardZunger1990,ZungerWeiFerreiraBernard1990}. It carries the name SQS, abbreviation of \emph{Special Quasirandom Structures}. The idea there is to sample random atomistic systems by selecting systems of finite size that exhibit all the suitable statistical properties usually only obtained in the asymptotic, bulk limit. For instance, one imposes on each and every finite size sample considered to have the correct volume fraction (proportion) of species (atoms, molecules, ions, etc), and/or the correct correlations present in the infinite limit. To some extent, the approach is actually a particular instance of a general approach, known in many different contexts, and popularized in political sciences because it is at the foundation of many techniques used for opinion polls: \emph{stratified sampling}. The approach has been adapted to the homogenization context in~\cite{sqs,minvielle-these} to which we refer the reader for a more detailed presentation. 

\paragraph{Mathematical setting.} 

Before explicitly formulating the selection mechanism employed in the method, let us motivate it. Consider a formal expansion of the coefficient $A$ present in~\eqref{eq:pb0-stoc} and~\eqref{eq:correcteur-random} in terms of a small parameter (even though, as we did in the previous section, we will soon forget about the smallness assumption). Assume $A$ reads as
\begin{equation}
\label{eq:A_form}
A_\eta(x, \omega) = C_0(x, \omega) + \eta \ \chi(x, \omega) \ C_1(x, \omega)
\end{equation}
for some presumably small scalar coefficient $\eta$, where we assume that $C_0$ and $C_1$ have nice properties. Since $\eta$ is small, $A_\eta$ is intuitively a perturbation of the matrix-valued field $C_0$.

At least formally (and all this has been mathematically justified in our works), we may expand all quantities of homogenization theory in powers of the small parameter~$\eta$. In particular, the approximations~$w^N_\eta$ and $A^{\star,N}_\eta$ of, respectively, the corrector~$w_\eta$ and the homogenized matrix~$A^\star_\eta$ on the truncated domain $Q_N$, can be expanded in powers of $\eta$:
\begin{eqnarray}
\label{eq:monde_reel_pre}
w^N_\eta(\cdot,\omega) &=& w^N_0(\cdot,\omega) + \eta u^N_1(\cdot,\omega) + \eta^2 u^N_2(\cdot,\omega) + o(\eta^2),
\\
\label{eq:monde_reel}
A^{\star,N}_\eta(\omega) &=& A^{\star,N}_0(\omega) + \eta A^{\star,N}_1(\omega) + \eta^2 A^{\star,N}_2(\omega) + o(\eta^2).
\end{eqnarray}
Inserting the expansion~\eqref{eq:monde_reel_pre} in the truncated corrector problem~\eqref{eq:correcteur-random-N}, one easily obtains a hierarchy of equations, the first three of which reading as
$$
\left\{
\begin{array}{ccc}
- \hbox{\rm div}\,C_0 (p + \nabla w^N_0) = 0& \quad \text{in } Q_N, &\quad \text{$w^N_0$ is $Q_N$-periodic},
\\ \noalign{\vskip 3pt}
- \hbox{\rm div}\,C_0 \nabla u^N_1 = \hbox{\rm div}\,\left[ \chi C_1 (p+\nabla w^N_0) \right]& \quad \text{in } Q_N, &\quad \text{$u^N_1$ is $Q_N$-periodic},
\\ \noalign{\vskip 3pt}
- \hbox{\rm div}\,C_0 \nabla u^N_2 = \hbox{\rm div}\,\left[ \chi C_1 \nabla u^N_1 \right]& \quad \text{in } Q_N, &\quad \text{$u^N_2$ is $Q_N$-periodic}.
\end{array}
\right.
$$
We next deduce from~\eqref{eq:Astar-random-N} and~\eqref{eq:monde_reel_pre} that~\eqref{eq:monde_reel} holds, with
\begin{eqnarray*}
A^{\star,N}_0(\omega) \, p &=& \frac{1}{|Q_N|} \int_{Q_N} C_0(\cdot,\omega) (p + \nabla w^N_0(\cdot,\omega)),
\\ \noalign{\vskip 3 pt}
A^{\star,N}_1(\omega) \, p &=& \frac{1}{|Q_N|} \int_{Q_N} \chi(\cdot,\omega) C_1(\cdot,\omega) (p+ \nabla w^N_0(\cdot,\omega)) + \frac{1}{|Q_N|} \int_{Q_N} C_0(\cdot,\omega) \nabla u^N_1(\cdot,\omega),
\\ \noalign{\vskip 3 pt}
A^{\star,N}_2(\omega) \, p &=& \frac{1}{|Q_N|} \int_{Q_N} \chi(\cdot,\omega) C_1(\cdot,\omega) \nabla u^N_1(\cdot,\omega) + \frac{1}{|Q_N|} \int_{Q_N} C_0(\cdot,\omega) \nabla u^N_2(\cdot,\omega).
\end{eqnarray*}

We may now introduce the conditions that we use to select particular configurations of the environment within $Q_N$. For finite fixed $N$, we say that a configuration $\omega$ satisfies the SQS conditions of order up to $k$ if, for any $0 \leq j \leq k$, the coefficient $A^{\star,N}_j(\omega)$ in the expansion~\eqref{eq:monde_reel} exactly matches the corresponding coefficient~$A^\star_j$ of the analogous expansion of the exact homogenized matrix coefficient $A^\star_\eta$. More explicitly, we speak about the SQS condition of 
\begin{itemize}
\item order 0 if $A^{\star,N}_0(\omega) = A^\star_0$, that is to say, for any $p \in \RR^d$,
\begin{equation}
\label{eq:cond0}
\frac{1}{|Q_N|}\int_{Q_N} C_0(x, \omega) (p + \nabla w^N_0(x, \omega)) dx = \EE\left[\int_Q C_0 (p + \nabla w_0)\right].
\end{equation}
\item order 1 if $A^{\star,N}_1(\omega) = A^\star_1$, that is to say, for any $p \in \RR^d$,
\begin{multline}
\label{eq:cond1}
\frac{1}{|Q_N|}\int_{Q_N} \Big( \chi(x, \omega)C_1(x, \omega) (p+ \nabla w^N_0(x, \omega)) + C_0(x, \omega) \nabla u^N_1(x, \omega) \Big) dx \\ 
= \EE\left[ \int_Q \chi C_1 (p+ \nabla w_0) + C_0 \nabla u_1\right].
\end{multline}
\item order 2 if $A^{\star,N}_2(\omega) = A^\star_2$, that is to say, for any $p \in \RR^d$,
\begin{multline}
\label{eq:cond2}
\frac{1}{|Q_N|} \int_{Q_N} \Big( \chi(x, \omega)C_1(x, \omega) \nabla u^N_1(x, \omega) +C_0(x, \omega) \nabla u^N_2 (x, \omega) \Big) dx \\ 
= \EE\left[\int_Q\chi C_1 \nabla u_1 + C_0 \nabla u_2\right].
\end{multline}
\end{itemize}

It is easily observed that, using such particular configurations that satisfy the SQS conditions of order up to $k$, we have, in the perturbative setting considered here, $A^{\star,N}_\eta(\omega)-A^\star_\eta = o(\eta^k)$. Taking the expectation over such configurations therefore formally provides a more accurate approximation of $A^\star_\eta$. Of course, the purpose is to apply the approach \emph{beyond} the perturbative setting. It can be expected (and it is indeed observed) that selecting the configurations using these conditions may reduce the variance. 

To make the computation of the right-hand sides of the above conditions practical (since in theory they can only be determined using an asymptotic limit, and are therefore as challenging to compute in practice as $A^\star$ itself), we restrict the generality of our setting. We assume that, in~\eqref{eq:A_form}, $C_0(x, \omega) = C_0$ is a deterministic, \textit{constant} matrix, $C_1(x, \omega) = C_1(x)$ is a deterministic, $\ZZ^d$-periodic matrix, and that $\displaystyle \chi(x, \omega) = \sum_{k \in \ZZ^d} X_k(\omega) 1_{k+Q}(x)$, where $X_k(\omega)$ are identically distributed, not necessarily independent, bounded random variables. For simplicity, we also assume here that $\EE\left[ X_0 \right] = 0$, and refer to~\cite{sqs,minvielle-these} for more general cases. After a tedious but not complicated calculation, we obtain that the two conditions~\eqref{eq:cond1}--\eqref{eq:cond2} rewrite as
\begin{eqnarray}
\label{eq:cond12_un}
\displaystyle \frac{1}{|Q_N|} \sum_{k \in \ZZ^d \cap Q_N} X_k(\omega) &=& 0,
\\ \noalign{\vskip 3pt}
\label{eq:cond12_deux}
\displaystyle \frac{1}{|Q_N|} \sum_{k, j \in Q_N \cap \ZZ^d} X_k(\omega) X_j(\omega) I_{k,j}^N &=& \sum_{k \in \ZZ^d} \EE[X_0 X_k] I_k^\infty,
\end{eqnarray}
respectively, where $\displaystyle I_k^\infty = \int_{k+Q} C_1 \nabla \phi_1$ and $\displaystyle I_{k,j}^N = \int_{Q+j} C_1(x) \nabla \phi_1^N(x-k)dx$. In these expressions, $\phi_1$ and $\phi_1^N$ solve $ -\hbox{\rm div}\,\left[ C_0 \nabla \phi_1 \right] = \hbox{\rm div}\,\left[ \mathbf{1}_Q C_1p\right] $ in $\RR^d$, and $ -\hbox{\rm div}\,\left[ C_0 \nabla \phi_1^N \right] = \hbox{\rm div}\,\left[ \mathbf{1}_{Q} C_1p \right] $ in $Q_N$ with periodic boundary conditions, respectively. The conditions~\eqref{eq:cond12_un}--\eqref{eq:cond12_deux} are called the SQS~1 and SQS~2 conditions. On the other hand, in the particular setting chosen, condition~\eqref{eq:cond0} (SQS~0, in some sense) is easily seen to be systematically satisfied when $N$ is an integer and the truncated approximation of~\eqref{eq:correcteur-random} that is chosen is the periodic approximation. 

\paragraph{Selection Monte Carlo sampling.}

The classical Monte Carlo sampling consists in successively generating a random configuration~$\omega_m$, solving the truncated corrector problem~\eqref{eq:correcteur-random-N} for that configuration, and computing the empirical mean~$\displaystyle {\cal I}^M_{MC} = \frac1M \sum_{m=1}^M A^\star_N(\omega_m)$ as an approximation for~$A^\star$.

In our selection Monte Carlo sampling, we systematically test whether the generated configuration satisfies the required SQS conditions, up to a certain tolerance, and reject it if it does not, \emph{before} solving the corrector problem~\eqref{eq:correcteur-random-N} for that configuration and letting it contribute to the empirical mean.

In full generality, the cost of Monte Carlo approaches is usually dominated by the cost of draws, and therefore selection algorithms are targeted to reject as few draws as possible. In sharp contrast (and this is a very important point which we cannot emphasize enough), in the present context where boundary value problems are to be solved repeatedly, the cost of draws for the configuration is negligible compared to the cost of the solution procedure for such boundary value problems. Likewise, evaluating the quantities present in~\eqref{eq:cond12_un}--\eqref{eq:cond12_deux} is inexpensive. Therefore, the purpose of the selection mechanism is to limit the number of boundary value problems to be solved, even though this comes at the (tiny) price of rejecting many configurations. We also note that, as for any selection procedure, our selection may introduce a bias (i.e. a modification of the systematic error in~\eqref{eq:error-decomposition}). In practice, we have observed that the gain in variance dominates the bias possibly introduced by the selection approach, and that the approach provides a better approximation of $A^\star$ than the standard Monte Carlo approach.

We have studied the approach theoretically in~\cite{sqs,minvielle-these}. We have shown that the estimator provided (at least in the simplest variant of our approach) converges to the homogenized coefficient $A^\star$ when the truncated domain progressively extends to the whole space. The efficiency of the approach is also theoretically demonstrated for some particular and simple situations (such as, again, the one-dimensional setting). A comprehensive experimental study of the approach has been completed.

\medskip

\begin{figure}[htbp]
\psfrag{Red = SQS 1st order, Blue = SQS 2nd order, Black- = MC, Black-- = exact}{}
\psfrag{meanvalue}{}
\psfrag{supercell size}{\footnotesize \hspace{-2mm} $N$}
\centering
\includegraphics[width=2.9in]{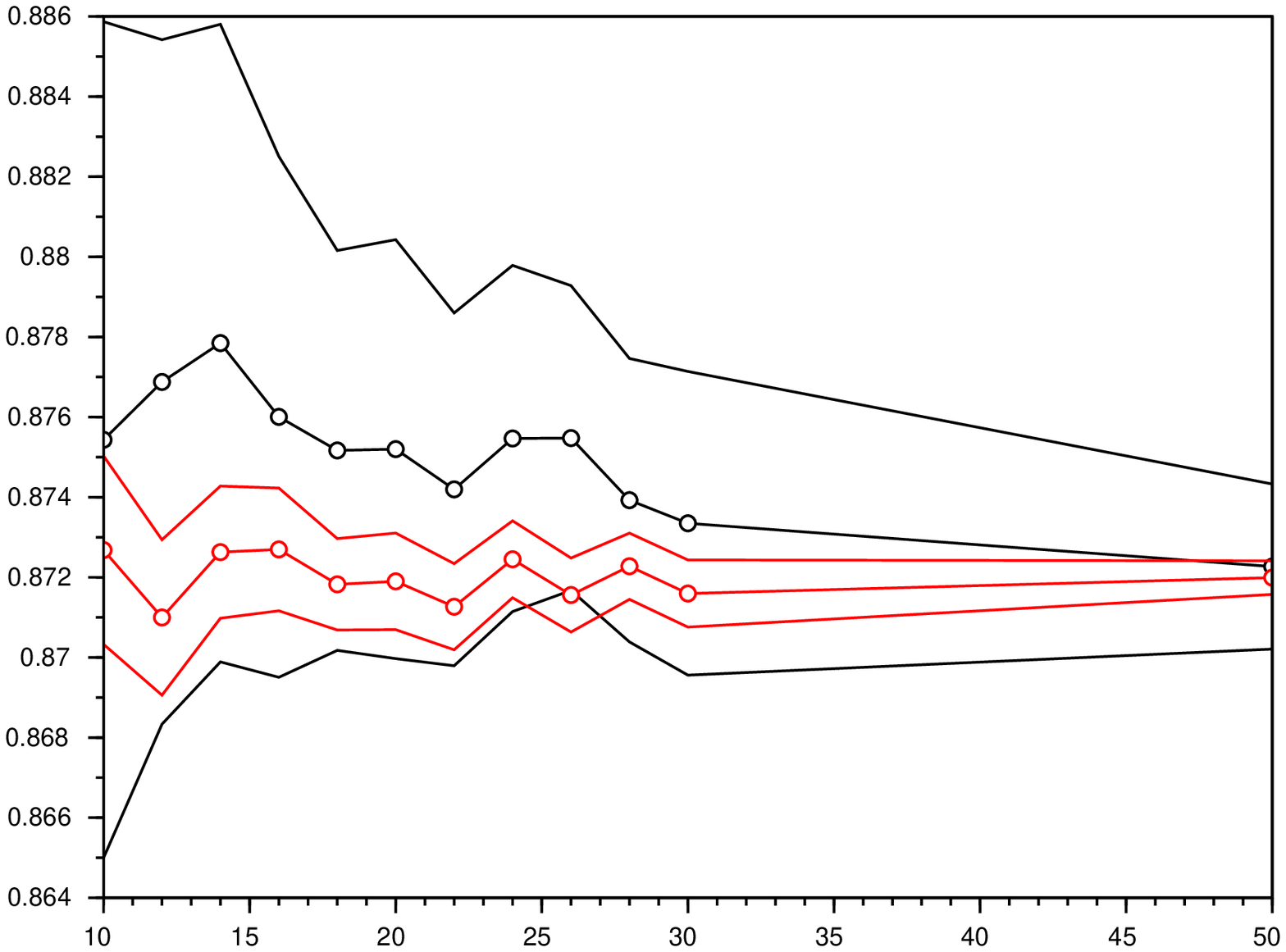}
\includegraphics[width=2.9in]{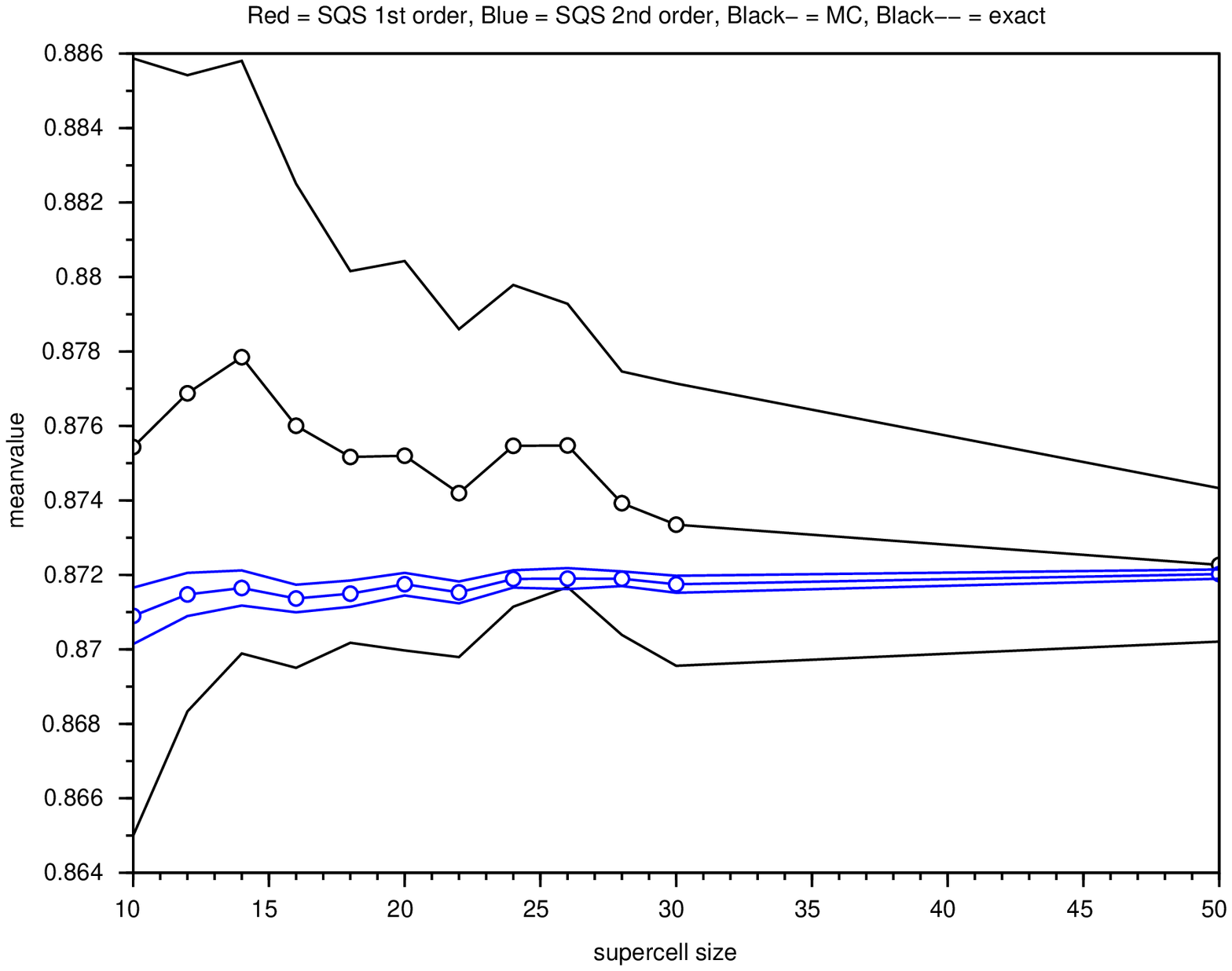}
\caption{Estimation of $A_{11}^\star$ together with its confidence interval (computed using $M=100$ i.i.d. realizations) as a function of $N$. Left: we use the classical MC approach (in black) and the SQS approach based on~\eqref{eq:cond12_un} (in red). Right: we use the classical MC approach (in black) and the SQS approach based on~\eqref{eq:cond12_un} and~\eqref{eq:cond12_deux} (in blue) (reproduced from~\cite{sqs}). \label{fig:SQS_un}}
\end{figure}

We include here a typical illustration of the efficiency of the approach. We again use a similar case as in Section~\ref{sssec:antithetic}. Considering only configurations that exactly satisfy~\eqref{eq:cond12_un}, we obtain the results shown on the left part of Figure~\ref{fig:SQS_un}. It is also possible, among the configurations that exactly satisfy~\eqref{eq:cond12_un}, to select configurations that satisfy as best as possible the condition~\eqref{eq:cond12_deux}. In practice, we generate 2000 configurations that exactly satisfy~\eqref{eq:cond12_un} and select among them the 100 configurations for which the difference between the left and the right-hand sides of~\eqref{eq:cond12_deux} is the smallest. We then obtain the results shown on the right part of Figure~\ref{fig:SQS_un}.

\medskip

In the case considered here (for which the contrast in the field $A$ is equal to 3), the variance is reduced by a factor 20 when using configurations that exactly satisfy~\eqref{eq:cond12_un}, and by a factor 300 if~\eqref{eq:cond12_deux} is enforced as well. To compare this variance reduction approach with the two previous ones, it is however needed to consider a case for which the contrast in $A$ is similar. In that case, the variance is reduced by a factor of 9 when using configurations that exactly satisfy~\eqref{eq:cond12_un}, and by a factor of 60 if~\eqref{eq:cond12_deux} is enforced as well.

\medskip

In all the test cases we have considered (see~\cite{sqs,minvielle-these} for details), we have observed that the systematic error is kept approximately constant by the approach (it might even be reduced), while the variance is reduced by several orders of magnitude. Such an efficiency is achieved at almost no additional cost with respect to the classical Monte Carlo algorithm.

\section{Robust Multiscale Finite Element approach for randomly perforated domains}
\label{sec:msfem}

\subsection{Presentation of our Crouzeix-Raviart type MsFEM approach with bubble functions}
\label{ssec:presentation2D}

\begin{figure}[htbp]
\psfrag{perf}{Perforations $B_\eps$}
\psfrag{bp}{Boundary $\partial B_\eps$ of the perforations}
\psfrag{dom}{Domain ${\mathcal D}_\eps$}
\centerline{
\includegraphics[width=7truecm]{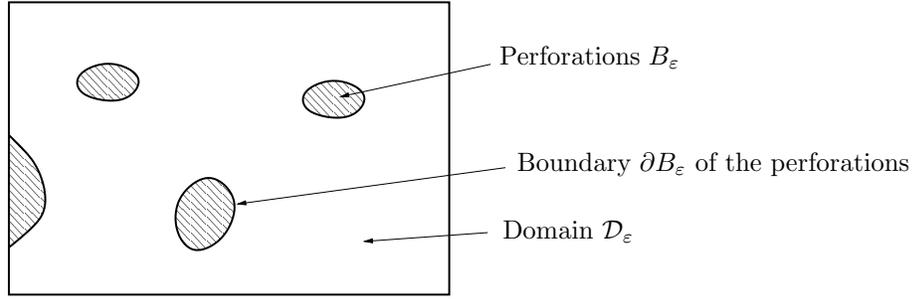}
}
\caption{The domain ${\mathcal D}$ contains (possibly random) perforations $B_\eps$. The perforated domain is ${\mathcal D}_\eps={\mathcal D} \setminus \overline{B_\eps}$. The boundary of ${\mathcal D}_\eps$ is the union of $\partial B_\eps \cap \overline{{\mathcal D}_\eps}$ (the part of the boundary of the perforations that is included in $\overline{{\mathcal D}_\eps}$) and of $\partial {\mathcal D} \cap \overline{{\mathcal D}_\eps}$. \label{fig:perforation}}
\end{figure}

We consider the problem~\eqref{eq:genP} in two dimensions, both for the analysis and for the numerical tests. We however emphasize that, of course, the approach can be applied to the three dimensional context and that, most likely, the theoretical analysis we review here can also be extended to the three dimensional case. Interestingly, our analysis in~\cite{loz1,loz3} on a similar problem was performed in both the two and three dimensional settings. We also assume, for simplicity, that the domain ${\mathcal D}$ is a polygonal domain. We define a mesh $\mathcal{T}_H$ on ${\mathcal D}$, i.e. a decomposition of ${\mathcal D}$ into polygons each of diameter at most $H$, and denote by $\mathcal{E}_H$ the set of all the internal edges of $\mathcal{T}_H$. Note that we mesh ${\mathcal D}$ and not the perforated domain ${\mathcal D}_\eps$. This allows us to use coarse elements (independently of the fine scale present in the geometry of ${\mathcal D}_\eps$), and leaves us with a lot of flexibility. The mesh does not have to be consistent with the perforations $B_\eps$. Some nodes may be in $B_\eps$, and likewise some edges may intersect $B_\eps$. See Figure~\ref{fig:perforation} for a representation of the perforated domain ${\mathcal D}_\eps$.

\paragraph{Functional space and variational approximation.}

To construct our MsFEM approach, we introduce the space
$$
W_H = \left\{
\begin{array}{c}
u \in L^2({\mathcal D}) \text{ such that } u|_T \in H^1(T)
\text{ for any } T\in \mathcal{T}_H,
\\ \noalign{\vskip 3pt}
\displaystyle
\int_E [[u]] =0 \text{ for all } E\in \mathcal{E}_H,
\quad u=0 \text{ in $B_\eps$ and on $\partial {\mathcal D}$}
\end{array}
\right\}
$$
where~$[[u]]$ denotes the jump of~$u$ across an edge, and the subspace
$$
W_H^0 = \left\{
\begin{array}{c}
\displaystyle
u \in W_H \text{ such that } 
\int_E u =0 \text{ for all } E\in \mathcal{E}_H 
\\ \noalign{\vskip 3pt}
\displaystyle
\text{ and } 
\int_T u =0 \text{ for all } T\in \mathcal{T}_H
\end{array}
\right\}
$$
of $W_H$. We define the Crouzeix-Raviart MsFEM space 
\begin{equation*}
V_H = \Big\{
u \in W_H \text{ such that } 
a_H(u,v) = 0 \text{ for all } v \in W_H^0 
\Big\}
\end{equation*}
as the orthogonal complement of $W_H^0$ in $W_H$, where by orthogonality we mean orthogonality for the scalar product $\dps a_H(u,v) = \sum_{T \in {\cal T}_H} \int_{T \cap {\mathcal D}_\eps} \nabla u \cdot \nabla v$. It is easy to see that we may write 
$$
V_H = \text{Span} \left\{ \Psi_T, \ \Phi_E, \quad E \in \mathcal{E}_H, \, T \in \mathcal{T}_H \right\},
$$ 
where the functions $\Psi_T$ and $\Phi_E$ are respectively defined as follows. 

First, for any mesh element $T$ that is not a subset of the perforations $B_\eps$ (i.e. $T \not\subset B_\eps$), we define $\Psi_T \in V_H$ as the solution to 
\begin{equation}
\label{eq:def_psi_bis}
-\Delta \Psi_T = 1 \ \text{in $T \setminus \overline{B_\eps}$},
\quad 
\Psi_T = 0 \ \text{in $T \cap B_\eps$},
\end{equation}
with, for each edge $\Gamma_i$ of $T$ which is an internal edge of $\mathcal{T}_H$, $\displaystyle \int_{\Gamma_i} \Psi_T = 0$ and $n \cdot \nabla \Psi_T = \lambda_i$ on $\Gamma_i$ for some constant $\lambda_i$. For any edge $\Gamma_i$ of $T$ which belongs to $\partial \Omega$, we set $\Psi_T = 0$ on $\Gamma_i$.

Second, for any internal edge $E$ that is not a subset of the perforations $B_\eps$, we denote by $T_E^1$ and $T_E^2$ the two triangles sharing this edge, set $T_E = T_E^1 \cup T_E^2$, and consider $\Phi_E \in V_H$ solution to
\begin{equation}
\label{eq:def_phi_bis}
-\Delta \Phi_E = 0 \ \text{in $T_E^1 \setminus \overline{B_\eps}$},
\quad 
-\Delta \Phi_E = 0 \ \text{in $T_E^2 \setminus \overline{B_\eps}$},
\quad 
\Phi_E = 0 \ \text{in $T_E \cap B_\eps$},
\end{equation}
with, for each edge $E' \subset \partial T_E$ which is an internal edge of $\mathcal{T}_H$, $\displaystyle \int_{E'} \Phi_E = 0$ and $n \cdot \nabla \Phi_E = \lambda_{E'}$ on $E'$ for some constant $\lambda_{E'}$ and $\displaystyle \int_E \Phi_E = 1$ and $n \cdot \nabla \Phi_E = \lambda_E$ on $E$ for some constant $\lambda_E$ (with an a priori different constant on the two sides of $E$). For any edge $E' \subset \partial \Omega$, we set $\Phi_E = 0$ on $E'$.

Third, for any mesh element $T \subset B_\eps$ (resp. any internal edge $E \subset B_\eps$), we set $\Psi_T \equiv 0$ (resp. $\Phi_E \equiv 0$).

\medskip

The MsFEM approximate solution is defined as the solution~$u_H \in V_H$ to 
\begin{equation*}
\forall v_H \in V_H, \quad a_H(u_H,v_H) = \int_{{\mathcal D}_\eps} f v_H.
\end{equation*}
 
\paragraph{Error estimate in the case of periodic perforations.}

Although this is not the focus of the present review, which specifically addresses approaches in the presence of randomness, we briefly mention for completeness that the approach we have just introduced enjoys suitable theoretical properties, at least in the periodic setting. One major reason why we mention this theoretical analysis is that, although there exist many works in the literature dealing with error estimates for MsFEM type approaches, none of them seems to apply (to the best of our knowledge) to problems set on perforated domains. In the periodic setting (and unfortunately we do not know how to extend our analysis to the random setting), we have established in~\cite{loz2} (in dimension~2, under technical conditions that are made precise in~\cite{loz2}, and using in particular arguments already present in~\cite{loz1,loz3}) the following error estimate:
\begin{equation}
|u-u_H|_{H^1_H({\mathcal D}_\eps)} \leq C \eps \, \left( \sqrt{\eps} + H + \sqrt{\frac{\eps}{H}} \right) \, \| f \|_{H^2({\mathcal D})},
\label{eq:mainresult}
\end{equation}
for some universal constant~$C$ independent from $H$, $\eps$ and $f$, but depending on the geometry of the mesh and other parameters of the problem. The norm used in the left-hand side of~\eqref{eq:mainresult} is the energy norm associated with the form $a_H$, namely $| v |_{H^1_H({\mathcal D}_\eps)} = \sqrt{a_H(v,v)}$ for any $v \in V_H + H^1_0(\Omega_\eps)$.

\medskip

Some important comments are in order. First of all, the size $\eps$ of the right-hand side of~\eqref{eq:mainresult} owes to the fact that, as it is well known theoretically (see~\cite{lions1980}), the size of the exact solution~$u$ (and thus that of the corresponding approximation~$u_H$) is~$\eps$ in $H^1$ norm. Taking this scaling factor into account, the actual rate of convergence for the numerical approach we design is therefore given by $\displaystyle \sqrt{\eps}+ H + \sqrt{\eps/H}$. Second of all, if the bubble functions $\Psi_T$ are not included in the discretization space, the estimate~\eqref{eq:mainresult} is replaced by
$$
|u-u_H|_{H^1_H({\mathcal D}_\eps)} \leq C \eps \, \left( \sqrt{\eps} + 1 + \sqrt{\frac{\eps}{H}} \right) \, \| f \|_{H^2({\mathcal D})}
$$
and the relative error does not converge to 0 in the limit when $\eps \to 0$ and next $H \to 0$, which is the relevant limit for our numerical experiments. This variant without bubbles is considered in our numerical tests for the sake of comparison.

The proof of~\eqref{eq:mainresult} is based upon a series of somewhat classical ingredients. We write
\begin{equation}
\label{eq:ineq_triangle}
u-u_H=(u-v_H)+(v_H-u_H),
\end{equation} 
where the function $v_H$ is defined as
\begin{equation*}
v_H(x) = \sum_{T \in {\cal T}_H} \Pi_H f \ \Psi_T(x) + \sum_{E \in {\cal E}_H} \left[ \int_E u \right] \ \Phi_E(x),
\end{equation*}
where the functions $\Psi_T$ and $\Phi_E$ are defined by~\eqref{eq:def_psi_bis} and~\eqref{eq:def_phi_bis} respectively and $\Pi_H f$ is defined as the $L^2$-orthogonal projection of $f$ on the space of piecewise constant functions. The most delicate contribution in the error to estimate is the first one in~\eqref{eq:ineq_triangle}. 

Given the standard finite element interpolation result
\begin{equation}
\label{eq:P1_EF}
\| f - \Pi_H f \|_{L^2({\mathcal D})}
\leq
C H \| \nabla f \|_{L^2({\mathcal D})}
\quad \text{for some $C$ independent of $H$ and $f$},
\end{equation}
and the fact that $v_H$ satisfies 
$$
\begin{array}{c}
\dps \text{for each triangle $T$}, \qquad -\Delta v_H = \Pi_H f \quad \text{in $T \setminus \overline{B_\eps}$},
\\ \noalign{\vskip 4pt}
\dps \text{for each internal edge $E$}, \qquad \int_E v_H = \int_E u,
\\ \noalign{\vskip 4pt}
v_H = 0 \quad \text{in $B_\eps$},
\end{array}
$$
we expect $v_H \in V_H$ to be close to the exact solution $u$. We establish that
\begin{multline}
| u-v_H |_{H^1_H({\mathcal D}_\eps)}^2 
=
\sum_{T \in \mathcal{T}_H} \int_{{\mathcal D}_\eps \cap T}
\nabla \left[ u - \eps^2 w(\cdot/\eps) f \right] \cdot \nabla \phi
+
\sum_{T\in \mathcal{T}_H} \int_{{\mathcal D}_\eps \cap T} 
\left[-\Delta (\eps^2 w(\cdot/\eps) f - v_H) \right] \phi
\\
+ \eps^2 \sum_{E \in \mathcal{E}_H}
\int_{E \cap {\mathcal D}_\eps} 
[[u-v_H]] \ n \cdot \nabla (w(\cdot/\eps) f),
\label{eq:decompo_error}
\end{multline}
where $\phi = u-v_H$ and $w$ denotes the periodic corrector associated to the periodic homogenization problem, that is the periodic solution to the problem $-\Delta w =1$ in $Y \setminus \overline{B}$, $w=0$ in the perforation $B$ (where $Y$ is the unit cell). We refer to~\cite{blp,Engquist-Souganidis,jikov} for more background on homogenization theory. 

We next successively bound the three terms of the right-hand side of~\eqref{eq:decompo_error}. Loosely speaking:
\begin{itemize}
\item the first term is small because of a classical homogenization result (see again~\cite{lions1980}) that states that $\eps^2 w(\cdot/\eps) f$ is indeed an accurate approximation of $u$,
\item the second term is small because, at the leading order term in $\eps$, the first factor in the integrand is equal to $-\Delta \left( \eps^2 w(\cdot/\eps) f \right) + \Delta v_H \approx f - \Pi_H f$ which is small due to~\eqref{eq:P1_EF}, 
\item the third term is estimated as an integral involving highly oscillatory periodic functions. 
\end{itemize}
This eventually yields the estimate~\eqref{eq:mainresult}. The details are made precise in~\cite{loz2}.

\subsection{Numerical tests}
\label{ssec:Numerical-tests}

We now solve~\eqref{eq:genP} for some particular settings, comparing our approach with other existing MsFEM type methods. Our primary focus is the possible intersections of the perforations with the edges of mesh elements. We address this question in Section~\ref{ssec:num-rob}, and check there the robustness of our approach with respect to the (temporarily deterministic, but still variable) location of the perforations: the fact that the mesh intersects, or does not intersect, the perforations has a very little influence on the (good) accuracy of our approach, in contrast to other approaches. Our numerical tests culminate with a couple of test-cases that show the excellent performance of our approach in the presence of randomness.

\medskip

For the sake of completeness, we begin our numerical tests with a comparison (for a given set of fixed perforations) with existing numerical MsFEM type approaches, in both cases when we enrich or not the variational space with bubble functions (in the case of our Crouzeix-Raviart MsFEM approach, the functions $\Psi_T$ defined by~\eqref{eq:def_psi_bis}). This is the purpose of Subsection~\ref{ssec:num-comp}.

\medskip

We mention that, in all our numerical experiments, we actually do not directly solve~\eqref{eq:genP} but a penalized version of this problem on ${\cal D}$ itself (see~\cite{loz2} for details). 

\medskip

We mesh ${\mathcal D}$ and not the perforated domain ${\mathcal D}_\eps$. Some nodes may thus be in the perforations $B_\eps$, some edges may intersect $B_\eps$, etc. We then face difficulties when using the MsFEM variant with linear boundary conditions (a variant that we compare here with our approach) to directly solve~\eqref{eq:genP}. Indeed, properly defining the MsFEM basis functions e.g. when edges intersect $B_\eps$ would not be straightforward. Similar difficulties arise for the oversampling variant. For this reason, we do not consider~\eqref{eq:genP}, but the penalized version of that problem as mentioned above. In contrast, note that the Crouzeix-Raviart variant we introduce here can be used either on~\eqref{eq:genP} or on its penalized version. 

\subsubsection{Comparison with existing approaches}
\label{ssec:num-comp}

We solve~\eqref{eq:genP} on the domain ${\mathcal D}=(0,1)^2$, with the right-hand side~$\displaystyle f(x,y) =\sin \frac{\pi x}{2} \, \sin \frac{\pi y}{2}$, and we take $B_\eps$ the set of discs of radius $0.35\eps$ periodically located on the regular grid of period $\eps =0.03$. The reference solution is computed on a mesh of size $1024 \times 1024$.

\medskip

The approaches we compare our approach with are the following four respective approaches:
\begin{itemize}
\item the standard Q1 finite element method on the coarse mesh of size $H$. Of course, we do not expect that method to perform well for this multiscale problem and we only consider it as a ``normalization''.
\item the MsFEM with linear boundary conditions. Although this method is now a bit outdated, it is still considered as the primary MsFEM approach, upon which all the other variants are built.
\item the MsFEM with oscillatory boundary conditions. This variant is restricted to the two-dimensional setting. It uses boundary conditions provided by the solution to the oscillatory ordinary differential equation obtained by taking the trace of the original equation on the edge considered. The approach performs fairly well on a number of cases, although it may also fail. 
\item the variant of MsFEM using oversampling. This variant is often considered as the ``gold standard'', although it includes a parameter (the oversampling ratio), the value of which should be carefully chosen. When this parameter is taken large, the method becomes (possibly prohibitively) expensive. 
\end{itemize}
In addition, we consider for each of those approaches, and for our specific Crouzeix-Raviart type approach, two variants: one with, and the other without a specific enrichment of the basis set elements using bubble functions. For all approaches but the Crouzeix-Raviart type approach that we propose, the bubble $\Psi$ on the quadrangle $Q$ is defined as the solution to
$$
-\Delta \Psi = 1 \text{ on $Q \cap {\mathcal D}_\eps$}, 
\quad \Psi = 0 \text{ on $\partial(Q \cap {\mathcal D}_\eps)$}.
$$
For our Crouzeix-Raviart approach, we recall that the bubble $\Psi$ is defined by~\eqref{eq:def_psi_bis}.

\medskip

For a given mesh size $H$, the cost for computing the basis functions (offline stage) varies from one MsFEM variant to the other. However, for a fixed $H$, all methods without (respectively, with) bubble functions essentially share the same cost to solve the macroscopic problem on ${\mathcal D}$ (online stage). More precisely, for a given cartesian mesh, and when using variants including the bubble functions, there are 1.5 times more degrees of freedom in our Crouzeix-Raviart approach than in the three alternative MsFEM approaches mentioned above. Since a logarithmic scaling is used for the x-axis in the figures below, this extra cost does not change the qualitative conclusions we draw below.

\medskip

\begin{figure}[htbp]
\centering
\includegraphics[width=7.4truecm]{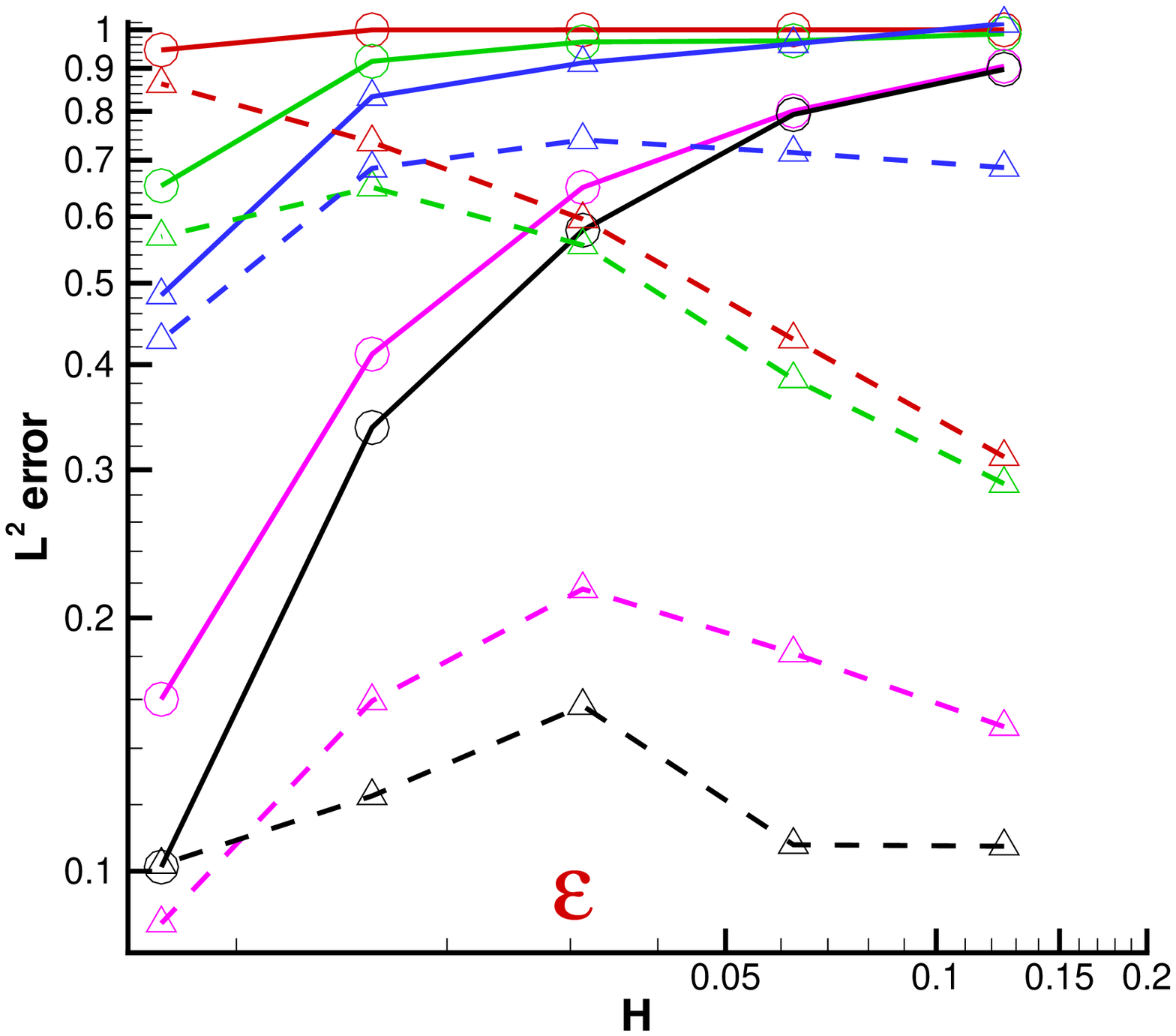}
\includegraphics[width=7.4truecm]{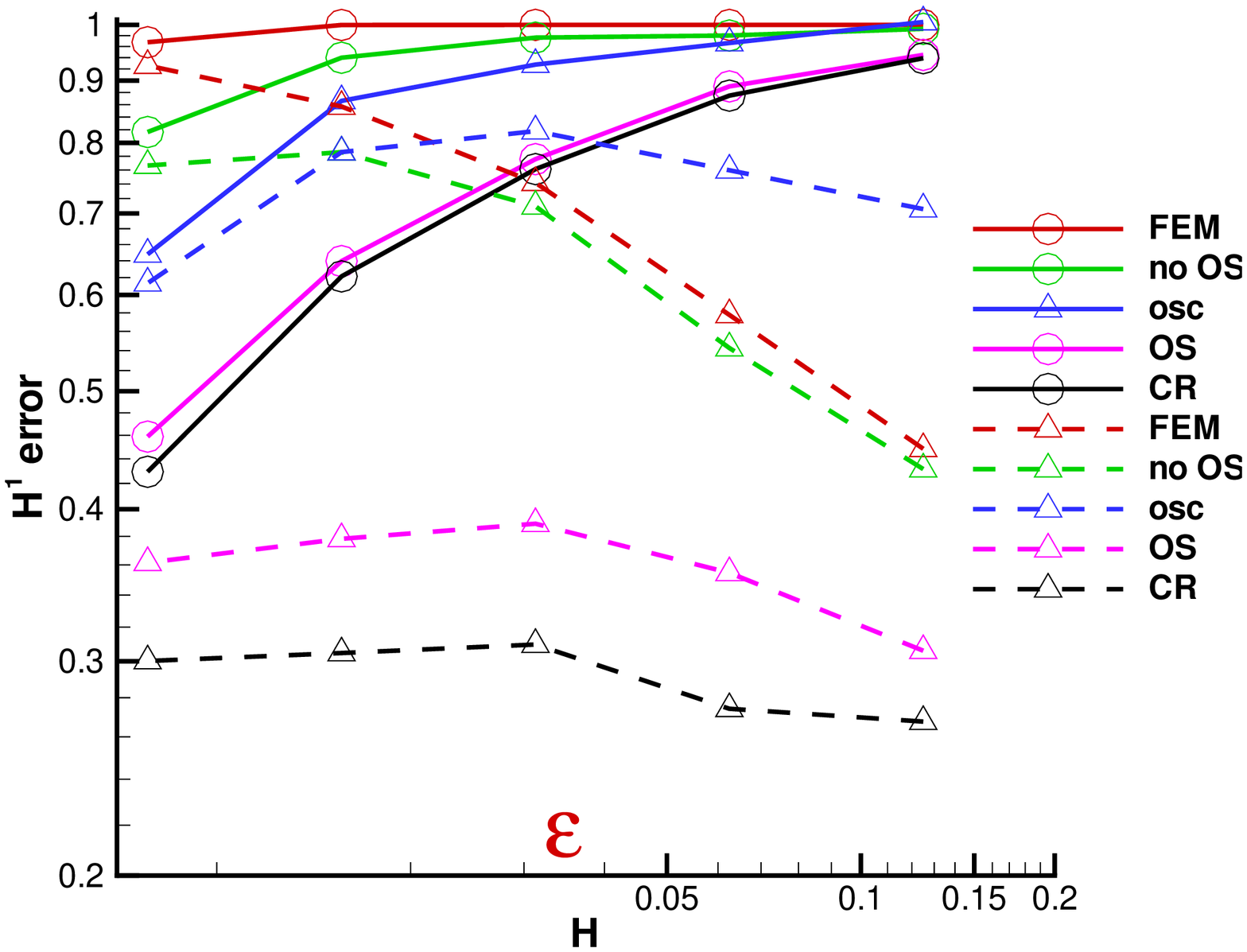}
\caption{Relative ($L^2$, left, and $H^1$-broken, right) errors with various approaches in the regimes $H \simeq \eps$ and $H \gtrsim \eps$: FEM -- the standard Q1 finite elements, no OS -- MsFEM with linear boundary conditions, osc -- MsFEM with oscillatory boundary conditions, OS -- MsFEM with oversampling (where the size of the quadrangles used to compute the basis functions is $3H \times 3H$), CR -- the Crouzeix-Raviart type MsFEM approach we propose. Results for all these methods are represented by solid lines. The dashed lines correspond to the variants of these methods where we enrich the finite element spaces using bubble functions (reproduced from~\cite{loz2}). \label{fig:errors}}
\end{figure}

The numerical results we have obtained in the regime where the meshsize $H$ is of the order of, or larger than, the parameter $\eps$ are presented on~Figure~\ref{fig:errors}. For all values of the meshsize~$H$, and for both the $L^2$ and the broken $H^1$ norms, a definite superiority of our approach over all other approaches is observed, and the interest of adding bubble functions to the basis set is, for each approach, also evident.

\medskip

A side remark is the following. On Figure~\ref{fig:errors}, we observe that, when using bubble functions, the error decreases as $H$ increases. This might seem counterintuitive at first sight. Note however that, when $H$ increases, the cost of computing each basis function increases, as we need to solve a local problem (discretized on a mesh of size $h$ controlled by the value of $\eps$) on a larger coarse element. In contrast to traditional FEM, increasing $H$ does not correspond to reducing the overall computational cost. For MsFEM approaches, increasing $H$ actually corresponds to decreasing the online cost but increasing the offline cost. The regime of interest is that of moderate values of $H$, for which the offline stage cost is acceptable. We only show the right part of Figure~\ref{fig:errors} (corresponding to large values of $H$, leading to a prohibitively expensive offline stage) for the sake of completeness.

\subsubsection{Robustness with respect to the location of the perforations}
\label{ssec:num-rob}

In this section, we perform a series of tests with the specific purpose of testing the robustness of the approach with respect to intersections of the mesh and the perforations. Given the flexibility of Crouzeix-Raviart type finite elements in terms of boundary conditions, we expect our approach to be particularly effective (and therefore considerably superior to other approaches) when some edges of the mesh happen to intersect perforations of the domain. The more frequent such intersections, the more important the difference. 

\paragraph{An extreme Best Case/Worst Case scenario.} 

We solve~\eqref{eq:genP} on the domain ${\mathcal D}=(0,1)^2$, with a constant right-hand side~$f=1$, and we take $B_\eps$ the set of discs of radius $0.2\eps$ periodically located on the regular grid of period $\eps =0.1$. We compute the reference solution, and consider 3 variants of MsFEM: the linear version, the oversampling version and the Crouzeix-Raviart version. The last three approaches are implemented in the variant that includes bubble functions in the basis set and applied with a mesh of size $H=0.2$.

\begin{figure}[htbp]
\centering
\includegraphics[width=7.4truecm]{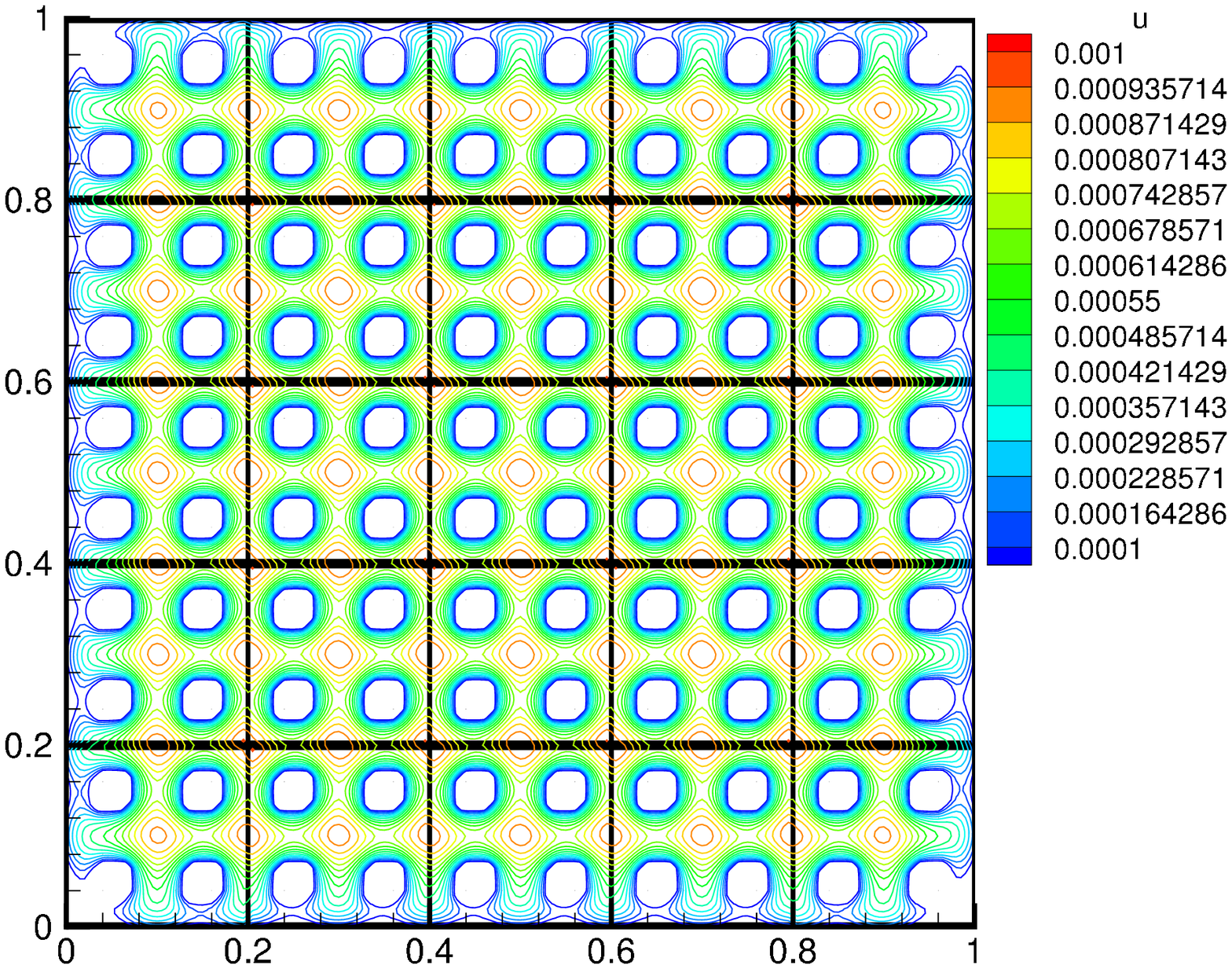}
\includegraphics[width=7.4truecm]{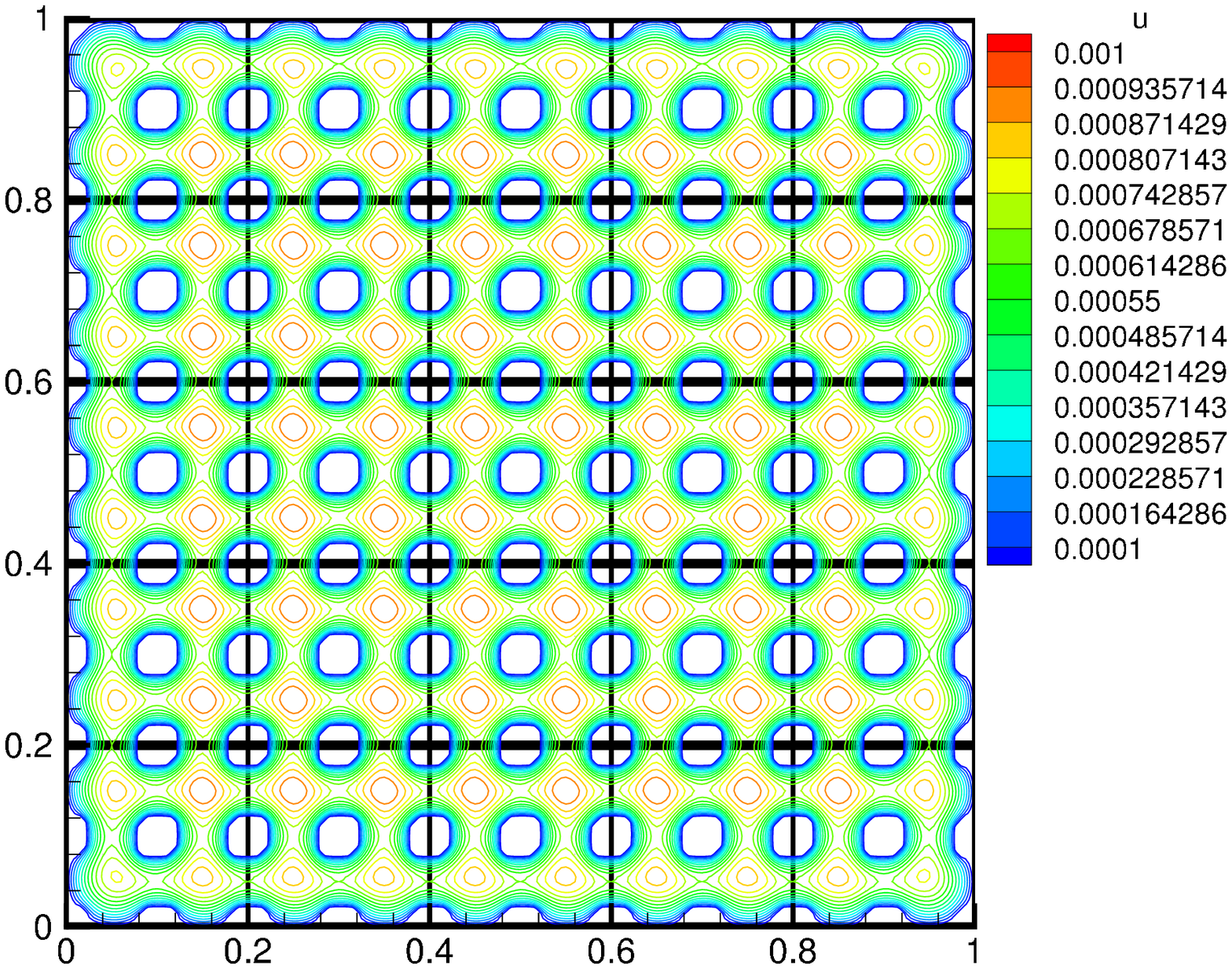}
\caption{Solution obtained with the proposed Crouzeix-Raviart type MsFEM (Left: Test 1; Right: Test 2) (reproduced from~\cite{loz2}). \label{fig:decalage-test1}}
\end{figure}

We perform two sets of numerical experiments. They are identical except for what concerns the relative position of the mesh with respect to the perforations. The difference between the two sets of tests is that, from one set of tests to the other one, the perforations are shifted by $\eps/2$ in the directions $x$ and $y$. In our Test~1, no edge intersects any perforation, while, on our Test~2, many edges actually intersect perforations (see Figure~\ref{fig:decalage-test1}). To some extent, the situation of Test 1 is the best case scenario (where as few edges as possible intersect the perforations) and the other situation is the worst case scenario. 

The numerical solutions are shown on Figure~\ref{fig:decalage-test1}. The numerical errors, computed both in $L^2$ and $H^1$-broken norms, are correspondingly displayed on Tables~\ref{table:test1} and~\ref{table:test2}, respectively. More than the actual values obtained for each case, this is the trend of difference between Table~\ref{table:test1} and Table~\ref{table:test2} that is the practically important feature. A comparison between the two tables indeed show that, qualitatively and in either of the norms used for measuring the error, the linear version and the oversampling version of MsFEM are both much more sensitive to edges intersecting perforations than the Crouzeix-Raviart version of MsFEM. In particular, the gain of our approach with respect to the linear version of MsFEM is much higher in our Test 2 (which is, from the geometrical viewpoint, the worst case scenario) than in Test 1. This confirms the intuition of a better flexibility of our approach. This also allows for expecting a much better behavior of that approach for non-periodic multiscale perforated problems for which it is extremely difficult, practically, to avoid repeated intersections of perforations with mesh edges. This is confirmed by the numerical experiments we perform in the next paragraph.

\medskip

\begin{table}[htbp]
\centering{
\begin{tabular}{l|c|c}
& $L^2$ error (\%) & $H^1$ error (\%) \\ \hline
MsFEM with linear conditions & 16 & 32 \\ \hline
MsFEM with oversampling & 20 & 38 \\ \hline
Crouzeix-Raviart type MsFEM & 9 & 24
\end{tabular}
}
\caption{Relative errors for Test 1. \label{table:test1}}
\end{table}

\medskip

\begin{table}[htbp]
\centering{
\begin{tabular}{l|c|c}
& $L^2$ error (\%) & $H^1$ error (\%) \\ \hline
MsFEM with linear conditions & 28 & 52 \\ \hline
MsFEM with oversampling & 12 & 31 \\ \hline
Crouzeix-Raviart type MsFEM & 9 & 27
\end{tabular}}
\caption{Relative errors for Test 2. \label{table:test2}}
\end{table}

\paragraph{Random perforations.}

\begin{figure}[htbp]
\centering
\includegraphics[width=6.4truecm]{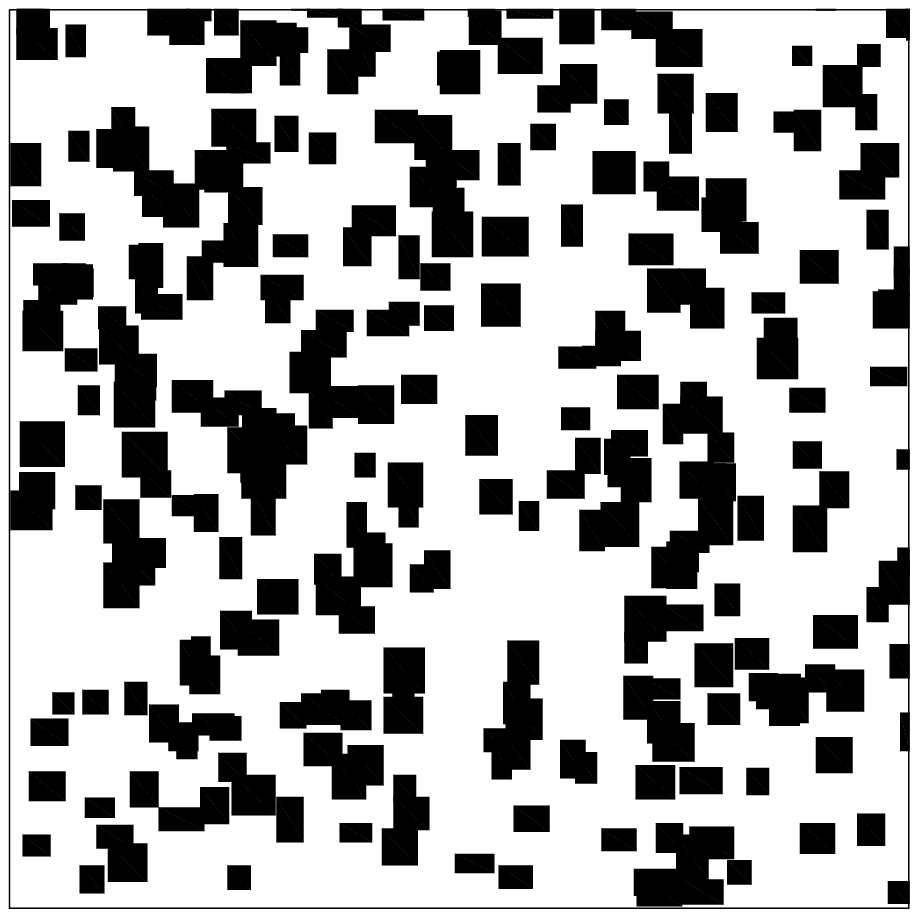}
\includegraphics[width=6.4truecm]{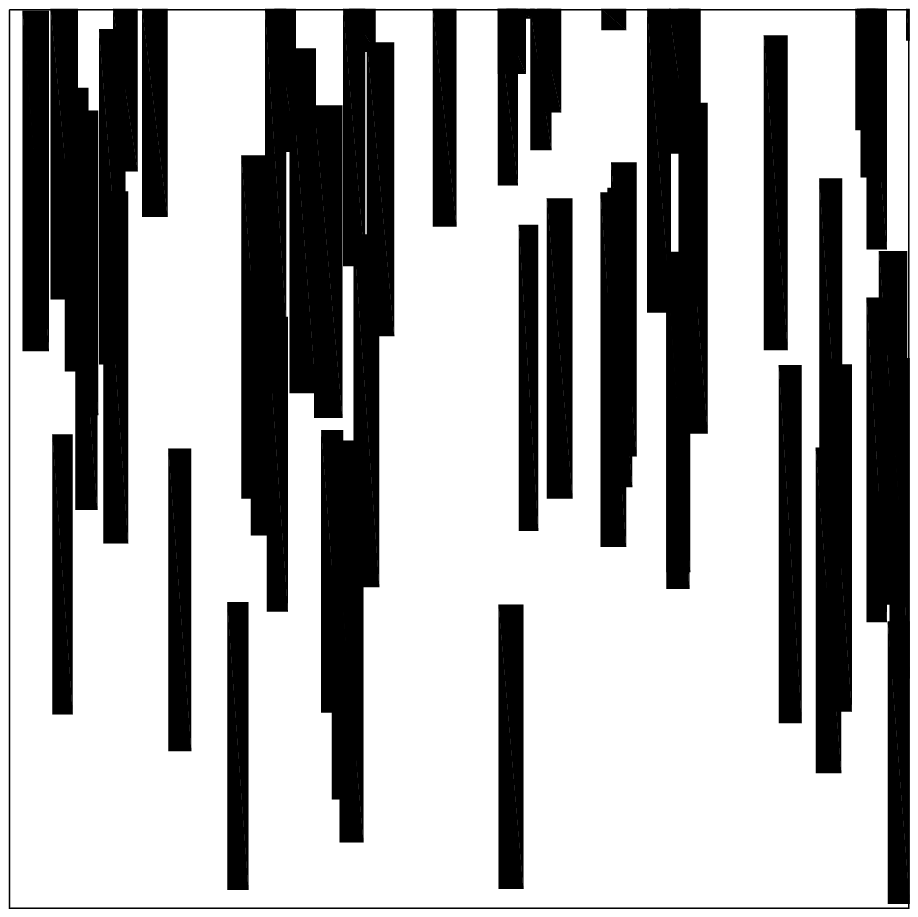}
\caption{Two examples of domains with random perforations (represented in black). Perforations have a rectangular shape, with a center randomly located in ${\mathcal D}=(0,1)^2$ according to the uniform distribution. Left: the perforations are 100 rectangles, the width and height of which are uniformly distributed between 0.02 and 0.05. Right: the perforations are 60 rectangles, the width (resp. the height) of which is uniformly distributed between 0.02 and 0.04 (resp. 0.02 and 0.4) (reproduced from~\cite{loz2}). \label{fig:ex-non-per}}
\end{figure}

We have tested several examples of randomly distributed perforations, two of them being shown on Figure~\ref{fig:ex-non-per}. For each of them, the domain ${\mathcal D} = (0,1)^2$ is meshed using quadrangles of size $H$, with $1/128 \leq H \leq 1/8$. The reference solution is again computed on a mesh of size $1024 \times 1024$.

\medskip

\begin{figure}[htbp]
\centering
\includegraphics[width=7.4truecm]{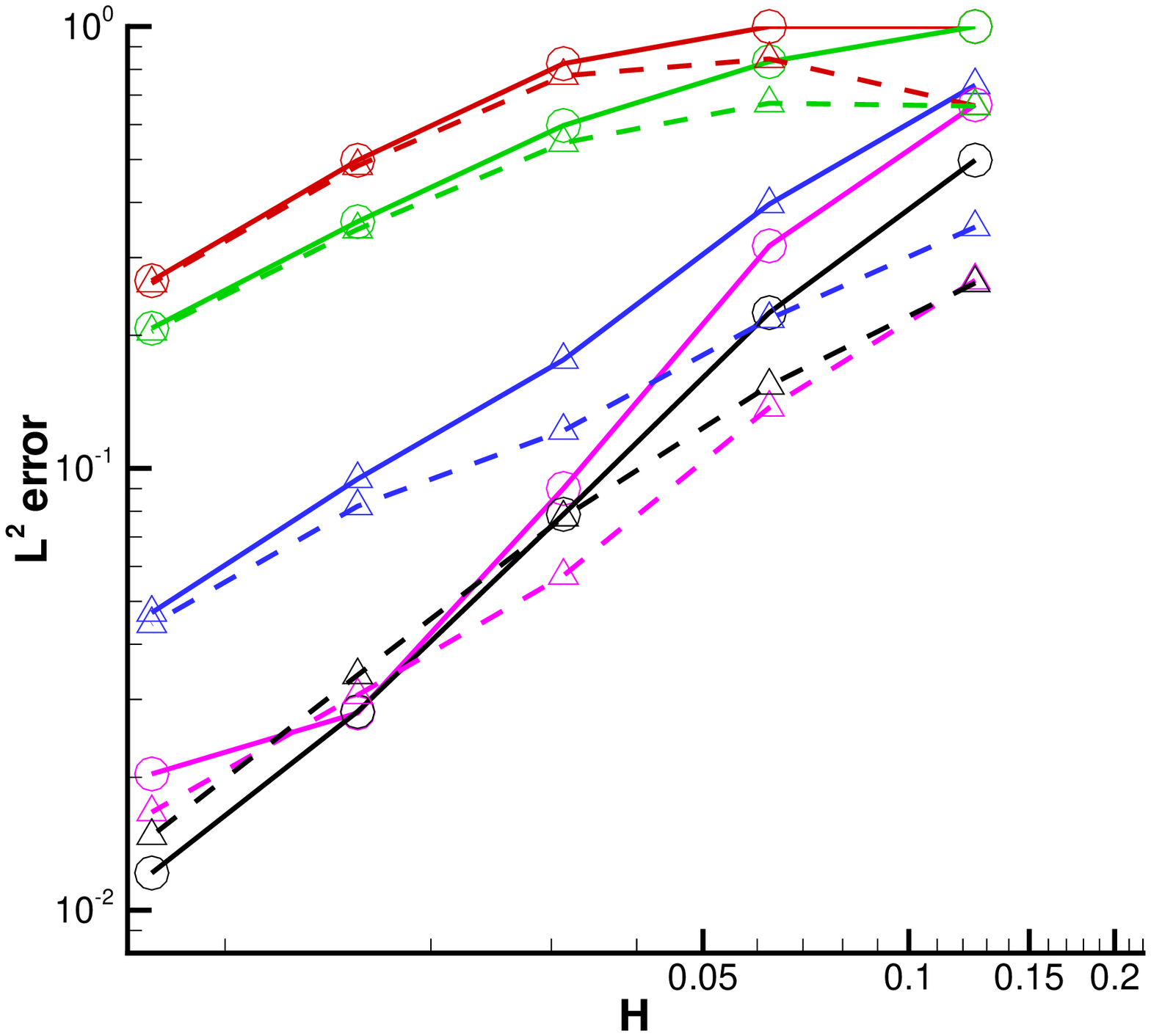}
\includegraphics[width=7.4truecm]{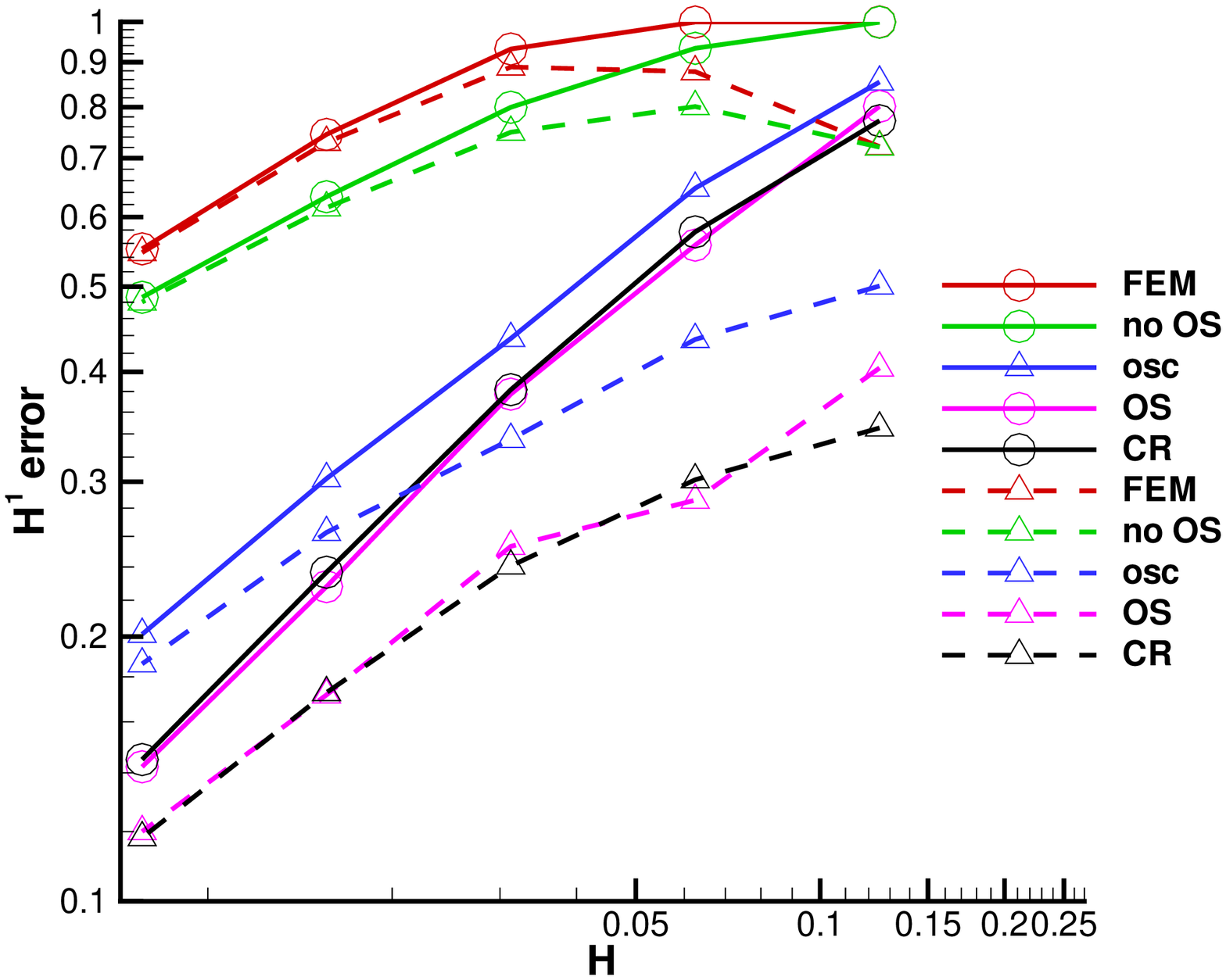}
\caption{Relative ($L^2$, left, and $H^1$-broken, right) errors with the same approaches as on Figure~\ref{fig:errors} for the test-case shown on the left part of Figure~\ref{fig:ex-non-per} (reproduced from~\cite{loz2}). \label{fig:errors-non-per3}}
\end{figure}

\begin{figure}[htbp]
\centering
\includegraphics[width=7.4truecm]{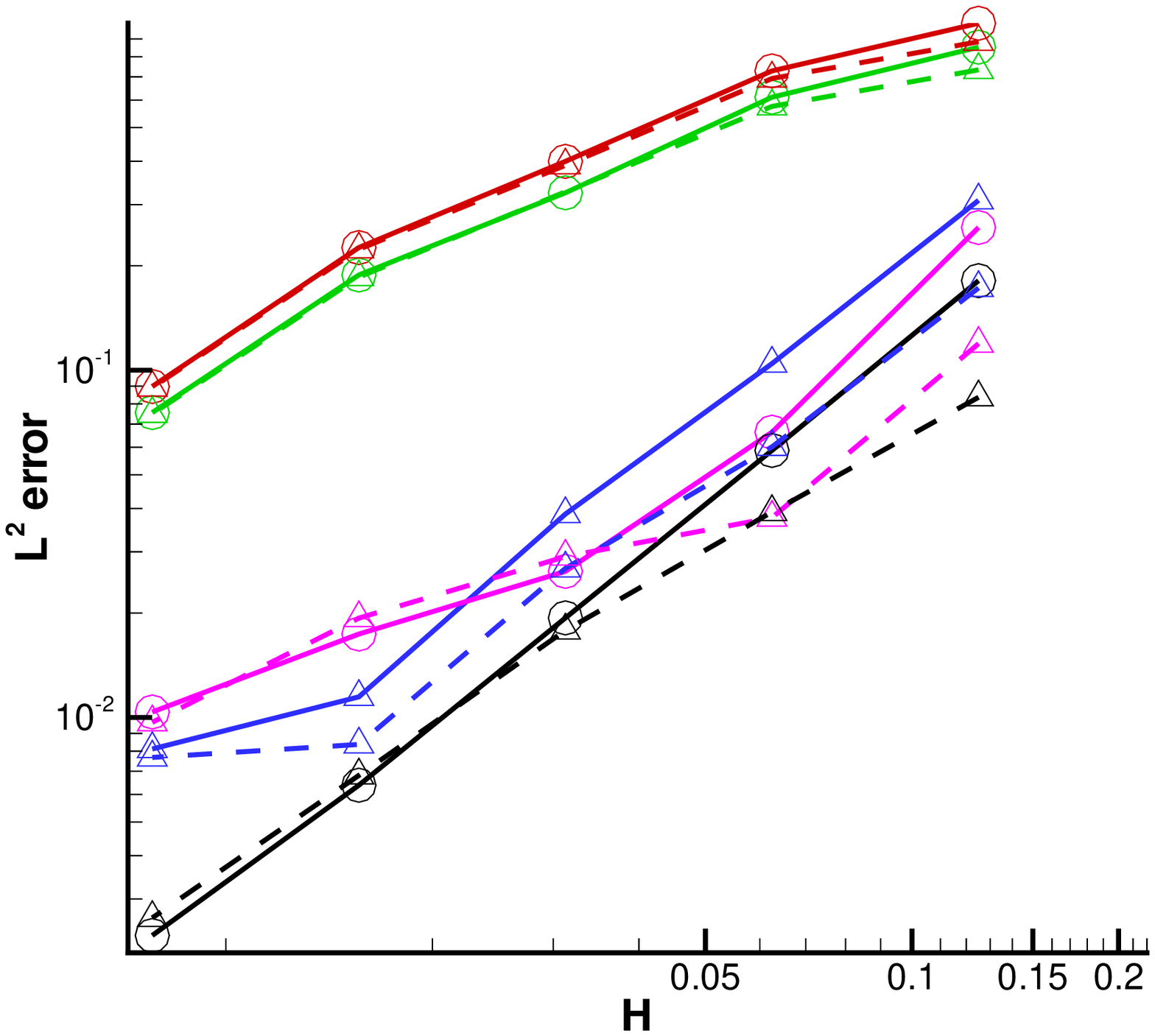}
\includegraphics[width=7.4truecm]{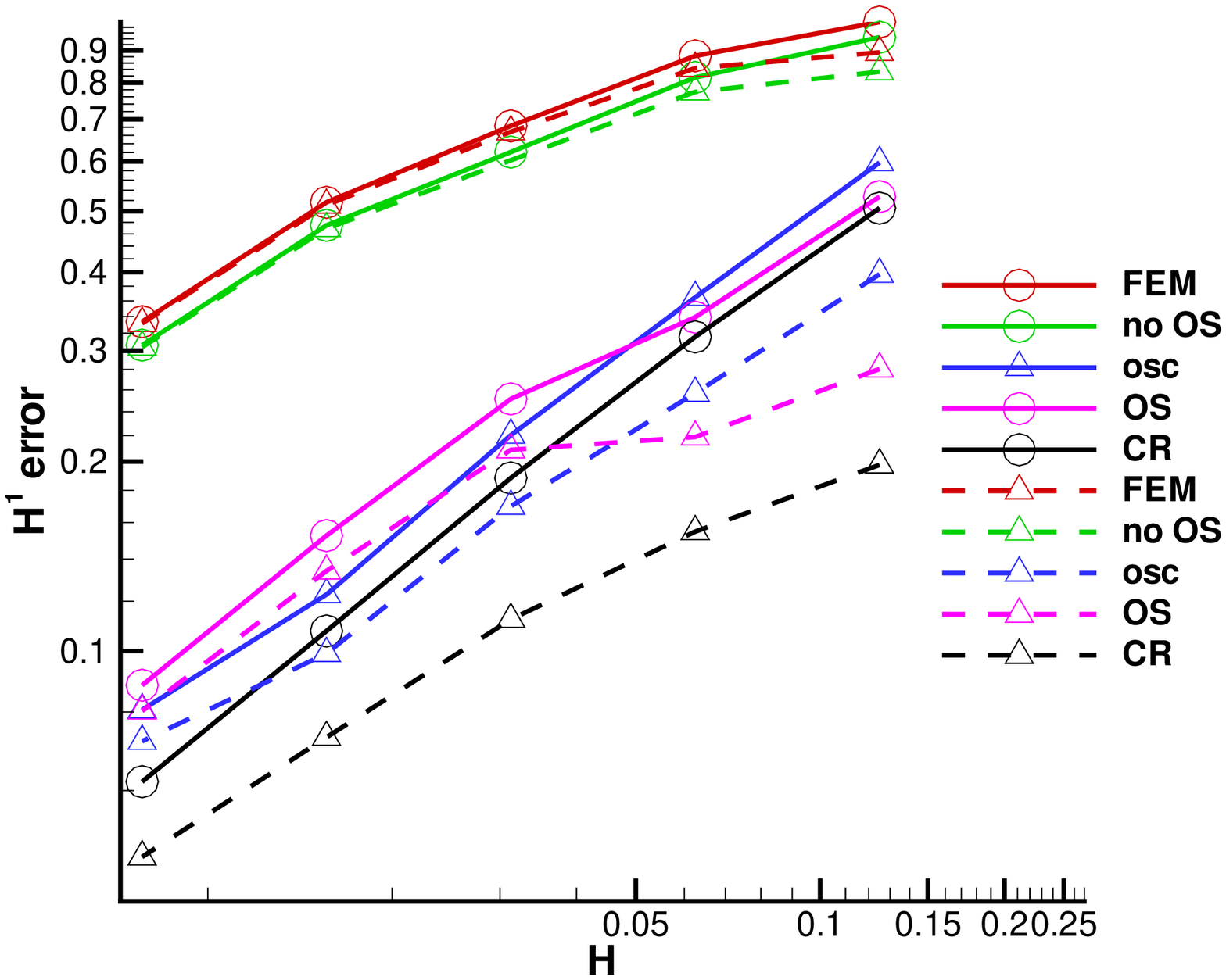}
\caption{Relative ($L^2$, left, and $H^1$-broken, right) errors with the same approaches as on Figure~\ref{fig:errors} for the test-case shown on the right part of Figure~\ref{fig:ex-non-per} (reproduced from~\cite{loz2}). \label{fig:errors-non-per5}}
\end{figure}

Errors are shown on Figure~\ref{fig:errors-non-per3} (resp. Figure~\ref{fig:errors-non-per5}) for the test-case shown on the left (resp. right) part of Figure~\ref{fig:ex-non-per}. Our approach provides results at least as accurate as, and often more accurate than the MsFEM approach with oversampling on quadrangles of size $3H \times 3H$. 

\medskip

\bigskip


\noindent{\bf Acknowledgements:} The work of the two authors is partially supported by EOARD under Grant FA8655-13-1-3061 and by ONR under Grant N00014-15-1-2777. The authors would like to thank their many friends and collaborators on the issues presented here, in particular X. Blanc (Paris 7), P.-L.~Lions (Coll\`ege de France), A.~Lozinski (Besan\c con) and W.~Minvielle. 

\bibliographystyle{plain}

\end{document}